\documentclass[11.5pt,reqno]{amsart}
\usepackage{amsthm, amsmath,amsfonts,amssymb,euscript,hyperref,graphics,color,slashed}
\usepackage{graphicx}
\usepackage{comment}
\usepackage{import}
\usepackage{tikz}
\usepackage{latexsym}
\usepackage{xcolor}

\def\ub{{\underline{u}}}
\def\Lb{\underline{L}}
\def\Cb{{\underline{C}}}
\def\Eb{\underline{E}}

\def\gslash{\mbox{$g \mkern -8mu /$ \!}}

\def\doubleint{\int\!\!\!\!\!\int}
\def\nablaslash{\mbox{$\nabla \mkern -13mu /$ \!}}

\definecolor{deepgreen}{cmyk}{1,0,1,0.5}

\newtheorem*{Main Theorem}{Main Theorem}
\newtheorem{theorem}{Theorem}[section]
\newtheorem{lemma}[theorem]{Lemma}
\newtheorem{proposition}[theorem]{Proposition}

\newtheorem{remark}[theorem]{Remark}

\numberwithin{equation}{section}

\begin{document}
\title[Large data regime for NLW]{On classical global solutions of nonlinear wave equations with large data}

\author[Shuang Miao]{Shuang Miao}
\thanks{Department of Mathematics, The University of Michigan,
Ann Arbor, MI 48109 U.S.A. shmiao@umich.edu}

\author[Long Pei]{Long Pei}
\thanks{Department of Mathematical Sciences, Norwegian University of Science and Technology. long.pei@math.ntnu.no}

\author[Pin Yu]{Pin Yu}
\thanks{Yau Mathematical Sciences Center, Tsinghua University, Beijing, China. pin@math.tsinghua.edu.cn}

\begin{abstract}
This paper studies the Cauchy problem for systems of semi-linear wave equations on $\mathbb{R}^{3+1}$ with 
nonlinear terms satisfying the null conditions. We construct future global-in-time classical solutions with arbitrarily 
large initial energy. The choice of the large Cauchy initial data is inspired by Christodoulou's characteristic initial data in his work \cite{Ch-08} on formation of black-holes.  The main innovation of the current work is that we discovered a relaxed energy ansatz which allows us to prove decay-in-time-estimate. Therefore, the new estimates can also be applied in studying the Cauchy problem for Einstein equations. 

\end{abstract}
\maketitle
\tableofcontents

\section{Introduction}
We consider the Cauchy problem of the following system of wave equations on $\mathbb{R}^{3+1}$:
\begin{equation}\label{Quadraticequation}
\Box\phi=Q(\nabla\phi,\nabla\phi).
\end{equation}
Here, $\Box =  -\partial_{t}^{2}+\triangle$ is the standard wave operator. The function $\phi$ is vector valued. In fact, $\phi$ stands for $N$ unknown functions $\phi^{I}$, $I=1,...,N$. The symbol $\nabla \phi$ denotes all possible $\partial_{\gamma} \phi^I$'s for $\gamma=0,1,2,3$ and
$I=1,2,\cdots,N$. The nonlinearity $Q(\nabla \phi, \nabla \phi)$ is a quadratic form in $\nabla\phi$ satisfying the \emph{null condition}, which will be specified later. The problem of constructing global-in-time solutions for small initial data has been studied intensively in the literature. The purpose of the current paper is to propose a large Cauchy data regime for \eqref{Quadraticequation} which also leads to global classical solutions.

\subsection{Historical results}
We discuss briefly the small data theory for \eqref{Quadraticequation} on $\mathbb{R}^{n+1}$. The idea is to use the decay mechanism of linear waves, i.e. solutions of $\Box\phi=0$, and treat the nonlinear problem as a perturbation of the linear waves. In dimensions greater than $3$, i.e. $n \geq 4$, the pointwise decay rate of linear waves is at least $t^{-3/2}$, which is integrable on for $t\geq 1$. This fast decay rate can be used to prove the small-data-global-existence results; see the pioneering works of Klainerman \cite{K-80} and \cite{K-84}. However, in $\mathbb{R}^{3+1}$, the pointwise decay of the linear wave is merely at the rate $t^{-1}$ which is not integrable. This weak decay rate is not enough to control the nonlinear interaction: F. John \cite{J-79} showed that there were quadratic forms (which do not satisfy the null condition) 
such that for arbitrarily small non-zero smooth data, solutions to \eqref{Quadraticequation} blow up in finite time.\\

The importance of the null condition was first observed in the breakthrough work \cite{K-85} by Klainerman, where he proved that small data lead to global-in-time classical solutions if the nonlinearity Q is a \emph{null form, which will be defined explicitly later}, or equivalently, satisfies the null condition. 
In \cite{Ch-86} Christodoulou obtained a similar result based on the conformal compactification of the Minkowski spacetime. Although the approaches in \cite{K-85} and \cite{Ch-86} are very different, both proofs rely on the cancellation structure of null condition, which is absent for general quadratic nonlinearities.\\ 

The idea of exploiting the cancellation structure of the null conditions can also be used to handle certain large data problems. In a recent breakthrough in general relativity, Christodoulou \cite{Ch-08} rigorously proved for the first time that black-holes can form dynamically from arbitrarily dispersed initial data. The key to this work was the new idea of the "short pulse method".  Roughly speaking, this is a choice of special large initial data, called \emph{short pulse} data, so that these large profiles can be propagated along the flow of Einstein vacuum equations. One of the key observations in the proof is still tightly related to the cancellation of the null structure: the profile is only large in certain components and these large components are always coupled with some small components so that their contributions are still manageable. Christodoulou's work has been generalized in \cite{K-R-09}  by Klainerman and Rodnianski. A key ingredient in their work is the relaxed propagation estimates which allows profiles with more large components.\\

The ideas used in \cite{Ch-08} and \cite{K-R-09} have been adapted to the main equation \eqref{Quadraticequation} by Wang and Yu to construct future-in-time global solutions with large initial data; see \cite{W-Y-12} and \cite{W-Y-13}. Their approach is indirect. The authors essentially impose the characteristic data on the past null infinity and solve the inverse scattering problem all the way up to a finite time to construct the initial Cauchy data. Very recently, Yang \cite{Yang-13} has obtained a global existence theorem for
semi-linear wave equations with large Cauchy initial energy. The largeness in \cite{Yang-13} is from a slower decay of the initial data at spatial infinity, but not from the short pulse method.\\

The aim of the current work is to study the global-in-time behavior of smooth solutions to \eqref{Quadraticequation} with short pulse data. We give short pulse Cauchy data directly (one should compare with the indirect approach of \cite{W-Y-12}) and prove that the data lead to future-global-in-time classical solutions for \eqref{Quadraticequation}. We remark that compared to the characteristic data approach in \cite{W-Y-12}, one of the main difficulties is to prove quantitative decay of the solution. This difficulty does not appear in \cite{W-Y-12}, because the data there are radiation fields given on the past null infinity, so that the decay rate is already explicitly given. We will give a more-detailed comparison of the present work and \cite{W-Y-12} after some necessary notations are introduced. 

\subsection{The short pulse data and main results}
We use $(x_0,x_1,x_2,x_3)$ to denote the standard Cartesian coordinates $(t,x,y,z)$ on $\mathbb{R}^{3+1}$. In particular, $\partial_0$ stands for $\partial_t$. Let $\phi: \mathbb{R}^{3+1} \rightarrow \mathbb{R}^N$ be a vector valued function, and we use $\phi^I$ to denote its components. We study the Cauchy problem for the following system of nonlinear wave equations
\begin{equation}\label{Main Equation}
\begin{split}
\Box\phi^{I}&=Q^I(\nabla \phi, \nabla \phi), \ \ \text{for } \ I=1,2,\cdots, N,\\
(\phi,\partial_t \phi)&\big|_{t=1} = (\phi_0, \phi_1).
\end{split}
\end{equation}
The 
nonlinearities $Q^I$ are \emph{null forms}, i.e. we can write $Q^I(\nabla \phi, \nabla \phi)$ as
\begin{equation*}
Q^I(\nabla \phi, \nabla \phi)=\sum_{0\leq\alpha,\beta\leq3, \atop 1\leq J,K\leq N}A^{\alpha\beta,\,I}_{JK}\partial_{\alpha}\phi^{J}\partial_{\beta}\phi^{K},
\end{equation*}
and for all \emph{null} vector $\xi \in \mathbb{R}^{3+1}$, i.e. $\xi=(\xi_0, \xi_{1}, \xi_{2}, \xi_3)$ satisfying $-\xi_{0}^{2}+\sum_{i=1}^3\xi^{2}_{i} =0$, the coefficient matrices $A^{\alpha\beta,\,I}_{JK}$ satisfy
\begin{align*}
\sum_{\alpha,\beta=0}^{3}A^{\alpha\beta,\,I}_{JK}\xi_{\alpha}\xi_{\beta}=0.
\end{align*}
For the sake of simplicity, instead of writing all the components, we shall always use $\phi$ to denote the $\phi^I$'s and use $Q(\nabla \phi, \nabla \phi)$ to denote the nonlinearity. In particular, we always write the main equation \eqref{Main Equation} as \eqref{Quadraticequation}. We remark that in order to simplify some of the expressions appearing in the proof of the main theorem, we give the initial data at $t = 1$ rather than $t = 0$. Because of the invariance of the equation under time translations, this is the same as giving data on $t=0$. 

\bigskip

Before describing the short pulse data, we introduce some notations: $r$ and $\theta$ are used to denote the usual radial and angular coordinates on $\mathbb{R}^{3}$. 
Let $\delta$ be a small positive constant which will be determined later. 
We identify the $t=1$ hypersurface with $\mathbb{R}^3$ and divide it into three parts: $$\big\{t=1\big\} = B_{1-2\delta} \cup \big(B_{1} -B_{1-2\delta}\big) \cup \big(\mathbb{R}^3-B_{1}\big),$$
where $B_r$ is the ball centered at the origin with radius $r$.

\bigskip

In the following, $f\lesssim g$ always means  there exists a constant $C$ such that $f\leq Cg$ holds.  
We consider the initial data $(\phi_0, \phi_1)$ on $\{t=1\}$ of \eqref{Main Equation} satisfying the following conditions:

\begin{itemize}
\item On $B_{1-2\delta}$, we set $(\phi_0, \phi_1) \equiv (0,0)$.
\item On $B_{1} -B_{1-2\delta}$,
\begin{equation}\label{data I}
\|\nabla^{k}\left(\phi_{1}+\partial_{r}\phi_{0}\right)\|_{L^{\infty}}\lesssim \delta^{1/2-k},
\end{equation}
and
\begin{equation}\label{data II}
\|\nabla^{k}\phi_{0}\|_{L^{\infty}}+\|\nabla^{k-1}\phi_{1}\|_{L^{\infty}}\lesssim \delta^{1/2-k}
\end{equation}
for any positive integer $1\leq k\leq 20$.
\item On $\mathbb{R}^3-B_{1}$, $(\phi_0, \phi_1) \equiv (0,0)$.
\end{itemize}
In particular, the following data satisfies \eqref{data I} and \eqref{data II}:
\begin{align}\label{data example}
\phi_{0}(r,\theta)=\delta^{1/2}\psi_{0}\left(\frac{1-r}{2\delta},\theta\right),\quad \phi_{1}(r,\theta)=\delta^{-1/2}\psi_{1}\left(\frac{1-r}{2\delta},\theta\right).
\end{align}
Here $\psi_{0}(s,\theta)$ and $\psi_{1}(s,\theta)$ are smooth functions supported in $(0,1)$ with respect to their first argument $s$. Moreover, $\psi_{0}$ and $\psi_{1}$ satisfy

\begin{align}\label{data constraint}
\|\left(\partial_{t}+\partial_{r}\right)\phi\|_{L^{\infty}(\Sigma_{1})}\lesssim\delta^{50},
\end{align}
and 

\begin{align}\label{data constraint 1}
\|\partial_{t}^{2}\phi-\partial_{r}^{2}\phi\|_{L^{\infty}(\Sigma_{1})}\lesssim \delta^{50}.
\end{align}
\eqref{data constraint} and \eqref{data constraint 1} can be achieved because we have the freedom to choose $\psi_{0}, \psi_{1}$ and there are two constraints to satisfy. 
Using an induction argument and in view of \eqref{data constraint}, \eqref{data constraint 1}, it is straightforward to see that for $k\leq 40$, we have

\begin{align}\label{L derivative estimate initial}
\|\left(\partial_{t}+\partial_{r}\right)^{k}\phi\|_{L^{\infty}(\Sigma_{1})}\lesssim\delta^{1/2}.
\end{align}

For a fixed value of $\theta$, the graph of $\phi_{0}$ versus the variable $r$ is as follows:

\includegraphics[width = 5.5 in]{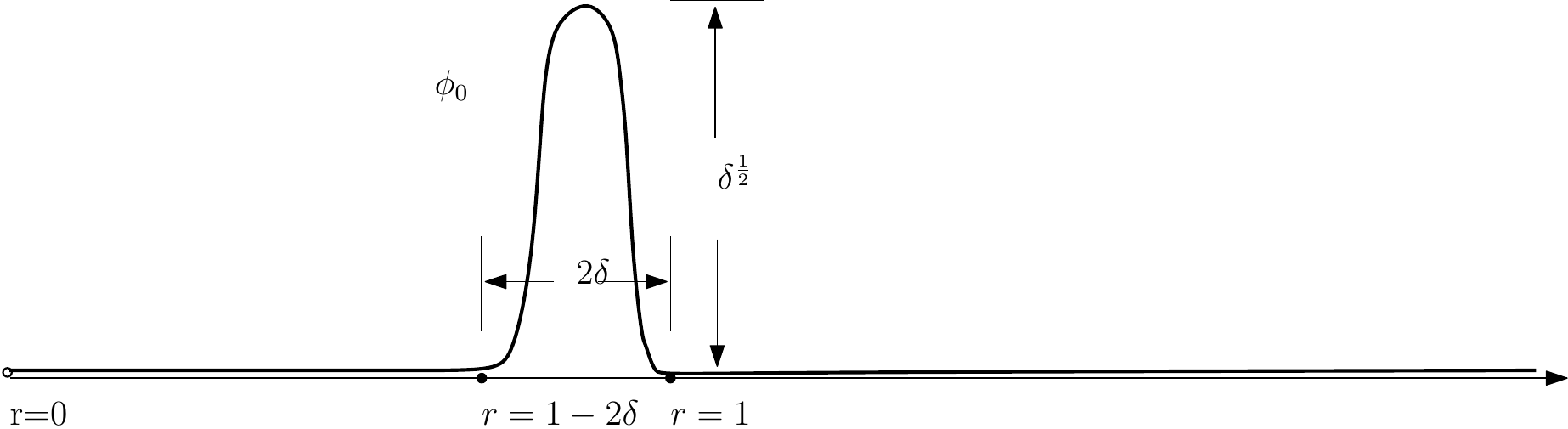}

The pulse-like shape of the graph explains the name
``short pulse" used for this data. 
The width of the pulse is $2\delta$ and its amplitude is $\delta^{\frac{1}{2}}$, which is very large relative to the width if $\delta$ is small.

\bigskip

The choice of 
$\phi_1(r,\theta)$ looks obscure and artificial in the above form. In fact, we have a natural geometric explanation of this choice, which can also serve as heuristics to understand why one expects a global-in-time solution.

\begin{remark}[Geometric / Physical interpretation]\label{geometric interpretation of the short pulse data}
In terms of the solution $\phi$, it is easy to observe from our initial data \eqref{data I} that
\begin{align}\label{small flux initial}
\left|(\partial_t + \partial_r) \phi |_{t=1}\right| \lesssim \delta^{1/2},\quad \left|\nablaslash\phi|_{t=1}\right|\lesssim \delta^{1/2}.
\end{align}
We will prove that if $\delta$ is small enough, the smallness indicated by \eqref{small flux initial} can be propagated by showing

\begin{align}\label{propagated small flux}
\|(\partial_{t}+\partial_{r})\phi\|_{L^{\infty}(\Sigma_{t})}\lesssim \delta^{1/2}t^{-2},\quad \|\nablaslash\phi\|_{L^{\infty}(\Sigma_{t})}\lesssim\delta^{1/4}t^{-2}.
\end{align}
Here we use $\Sigma_{t}$ to denote the hypersurface $\{t=t\}$. Recall that $L = \partial_t + \partial_r$ is the normal (with respect to the Minkowski metric!) of the outgoing light cones $t-r = \text{constant}$ in $\mathbb{R}^{3+1}$. If we integrate $|(\partial_t + \partial_r) \phi|^2+|\nablaslash\phi|^{2}$ on such an outgoing light cone $C$, the quantity
$$\int_{C}|(\partial_t + \partial_r) \phi|^2+|\nablaslash\phi|^{2}d\mu_C$$
measures the incoming energy through this light cone. Therefore, since $\delta$ will be eventually very small, the choice of $\phi_1$ is to keep the incoming energy as small as possible. Intuitively, we expect all the energy will be emanated in the outgoing direction so that the solution $\phi$ disperses.
\end{remark}

We now explain in what sense the short pulse data are large. It appears that the short pulse data is at least small in the $L^\infty$ sense due to the presence of the factor $\delta$. First of all, we notice that the $L^\infty$ norm is irrelevant since we may always add a constant to get a new solution for \eqref{Main Equation}. The size of the data should be measured at least on the the level of first derivatives. Secondly, we notice that, if we take derivatives in the $\partial_r$ direction many times, the data can be extremely large in the $L^\infty$ sense, because each $\partial_r$ derivative will bring out a $\delta^{-1}$ factor from the first argument of $\phi_0$ or $\phi_1$.

A more natural way to see the largeness of the data is to consider the energy spaces, i.e. the Sobolev spaces $H^k(\mathbb{R}^3)$. The critical $H^s$-exponent (with respect to scaling) of \eqref{Main Equation} is $\frac{3}{2}$. Therefore, the $0^{\text{th}}$ order energy $\mathcal{E}_0 = \int_{\mathbb{R}^3}|\nabla_{x}\phi_{0}|^{2}+|\phi_{1}|^{2}dx$ is subcritical and the $1^{\text{st}}$ order energy $\mathcal{E}_1=\sum_{i=1}^3\int_{\mathbb{R}^3}|\nabla_{x} \partial_i \phi_{0}|^{2}+|\partial_i \phi_{1}|^{2}dx$ is supercritical. Here we use $\nabla_{x}$ to denote spatial gradient.

\begin{remark}[Largeness of short pulse data]
We can compute the $0^{\text{th}}$ order energy $\mathcal{E}_0$ and the $1^{\text{st}}$ order energy $\mathcal{E}_1$ as follows:
\begin{align*}
\mathcal{E}_0 &\sim \|\nabla_{x} \phi_0\|^2_{L^2} +  \| \phi_1\|^2_{L^2} \sim 1,\\
\mathcal{E}_1 &\sim \|\nabla_{x}^2 \phi_0\|^2_{L^2} +  \|\nabla_{x} \phi_1\|^2_{L^2} \sim \delta^{-2}.
\end{align*}
Since $\mathcal{E}_0$ and $\mathcal{E}_1$ are subcritical and supercritical, respectively, we can not make both of them small by the scaling invariance of the equation. It is in this sense that the data are large at the level of energy.
\end{remark}

Moreover, for all $k\geq 0$, we can show that
\begin{equation*}
\mathcal{E}_k  =\int_{\mathbb{R}^3}|\nabla^{k+1}_{x} \phi_{0}|^{2}+|\nabla_{x}^k \phi_{1}|^{2}dx \sim \delta^{-2k}.
\end{equation*}
We note in passing that the higher order energies can be extremely large. Also, we remark that the symbol $\sim$ depends only on an absolute constant. 

\bigskip

We are now ready to state the main theorem of the paper:

\begin{Main Theorem} For any given pair of short pulse data $(\phi_0,\phi_1)$ as above, let us consider the following system of wave equations
\begin{align*}
&\Box\phi^{I}=Q^I(\nabla \phi, \nabla \phi), \ \ \text{for } \ I=1,2,\cdots, N,\\
&(\phi,\partial_t \phi)\big|_{t=1} = (\phi_0, \phi_1).
\end{align*}
where the $Q^I$'s are null forms.

Then there exists an absolute positive number $\delta_0$, so that for all $\delta < \delta_0$, the above Cauchy problem admits a unique smooth solution $\phi$ with lifespan $[1,+\infty)$. Moreover, when $t\rightarrow \infty$, the nonlinear wave $\phi$ scatters.
\end{Main Theorem}

\subsection{Notations}
We review the basic geometry of Minkowski space $\mathbb{R}^{3+1}$. In particular, we discuss the standard double null (cone) foliations on $\mathbb{R}^{3+1}$ which will play a central role for the energy estimates.

\bigskip

Let $r = \sqrt{x_1^2 + x_2^2 +x_3^2}$. We define two optical functions $u$ and $\ub$ as follows
\begin{equation*}
u = \frac{1}{2}(t-r),  \quad \ub = \frac{1}{2}(t+r).
\end{equation*}

For a given constant $c$, we use $C_{c}$ to denote the level surface $u=c$ with an extra constraint that $t \geq 1$ (since we will construct a future-global-in-time solution starting from the initial hypersurface $\{t=1\}$). According to the different value of $u$, we use also $C_u$ to denote these hypersurfaces. These are called outgoing light cones. Thus, $\{C_u\big|u\in \mathbb{R}\}$ defines a foliation of $\mathbb{R}^{3+1}_{\, t\geq 1}$. We also call this foliation null because each leaf $C_u$ is a null hypersurface with respect to the Minkowski metric.

Similarly, using the level sets of the optical function $\ub$, we define another null foliation of  $\mathbb{R}^{3+1}_{\, t\geq 1}$, denoted by $\{\Cb_{\ub}\big|\ub \in \mathbb{R}\}$. Each $\Cb_{\ub}$ is a truncated incoming light cone. The intersection $C_u \cap \Cb_{\ub}$ is a round $2$-sphere with radius $\ub - u$, denoted by $S_{\ub,u}$. We say that the two foliations $\{\Cb_{\ub}\big|\ub \in \mathbb{R}\}$ and $\{C_u\big|u\in \mathbb{R}\}$ form a \emph{double null foliation} of  $\mathbb{R}^{3+1}_{\, t\geq 1}$.

We recall that the following two null vector fields,
\begin{equation*}
 L = \partial_t + \partial_r,  \ \ \text{and} \ \ \Lb = \partial_t - \partial_r
\end{equation*}
are the normals of (also parallel to) $C_{u}$ and $\Cb_{\ub}$ respectively. In the following, null  pair always refers to the pair of two null vector fields $(L,\Lb)$.

The following picture depicts the outgoing null foliation $C_{u}$ of $\mathbb{R}^{3+1}_{\, t\geq 1}$:

\includegraphics[width = 5 in]{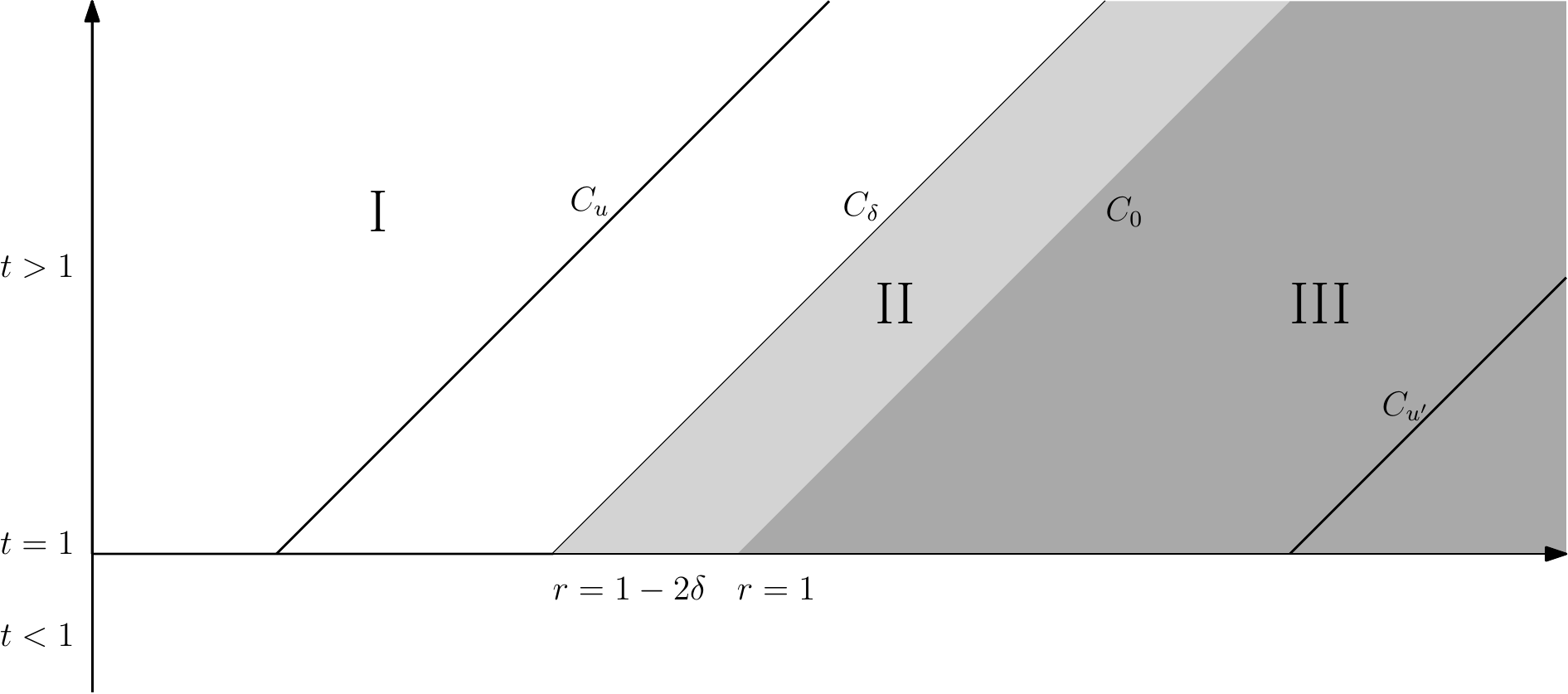}

Since the foliation is spherically symmetric, 
we only draw the $t$ and $r$ components in the schematic diagram. The other pictures in the paper should also be understood in this way. In the above picture, a $45^\circ$ line denotes an outgoing cone $C_u$.  Two outgoing cones $C_0$ and $C_\delta$ divide $\mathbb{R}^{3+1}$ into three regions: the small data region, i.e. region I in the picture, the short pulse region, i.e. the region with light grey color, region II in the picture, and the region III in the picture, i.e. the region with dark grey color.

\begin{remark}[Vanishing Property on 
$C_{0}$]
Recall that the short pulse data prescribed  on $\{t=1\}$ in the last subsection are identically zero for $r\geq 1$,
 therefore, according to the weak Huygen's principle, the solution of the main equation \eqref{Main Equation} vanishes identically in the region III (dark grey). In particular, the solution $\phi$ (if it exists) and its derivatives vanish on 
 $C_{0}$.
\end{remark}

We now pay more attention to the short pulse region (region II with light grey color). We use $D_{\ub,u}$ to denote the interior of the spacetime region enclosed by the hypersurfaces $\{t=1\}$, $C_{0}$, $C_{u}$ and $\Cb_{\ub}$, where $u \in [0,\delta]$ and $\ub \geq 1-u$. The following picture is a schematic diagram for all the notations introduced in the this section for the short pulse region.

\includegraphics[width = 4.5 in]{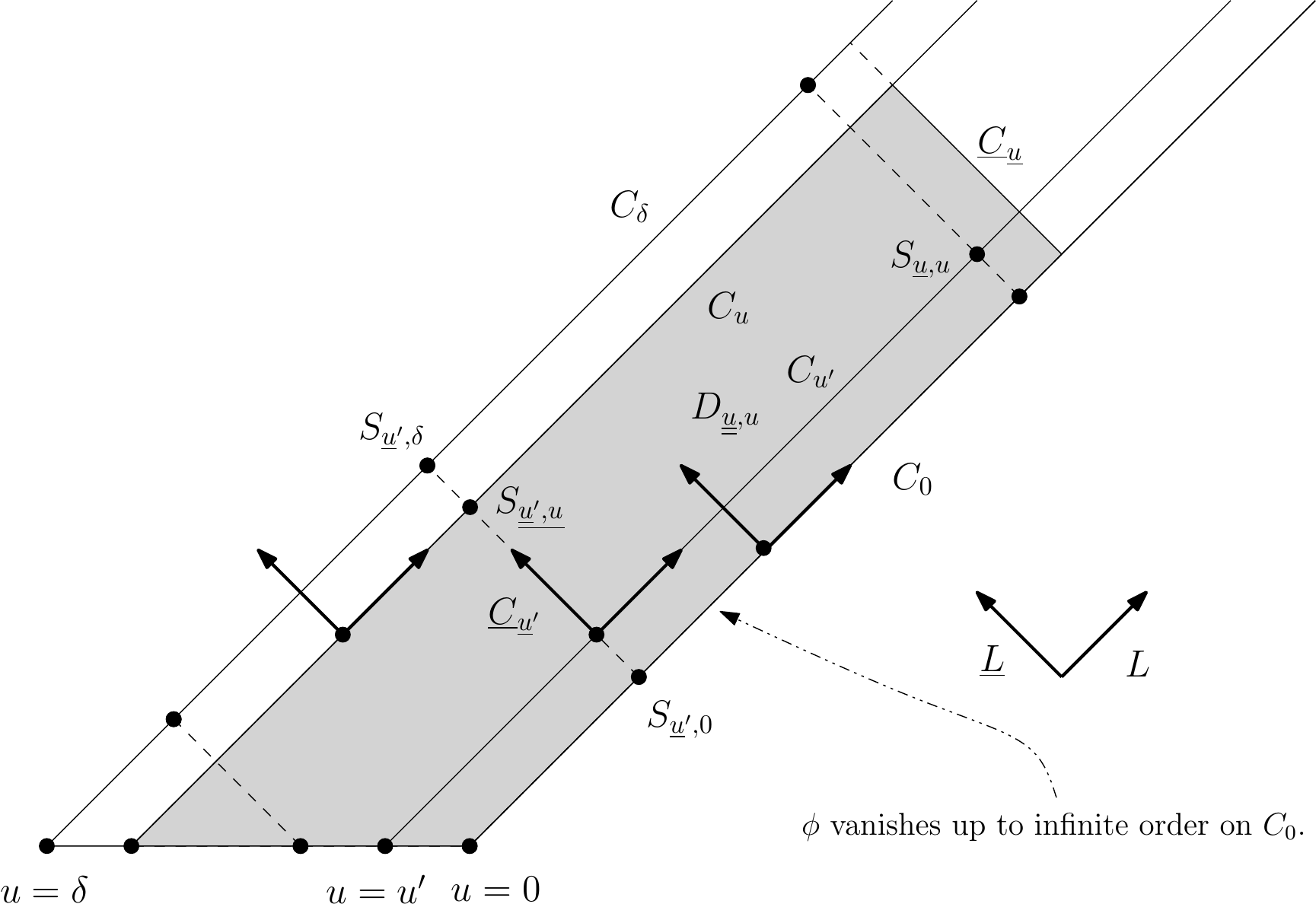}

A dashed $45^\circ$ segment  denote an incoming cone $\Cb_{\ub}$. A thickened black point denotes a $2$-sphere $S_{\ub,u}$. An orthogonal pair of arrow denotes the null vector pair $(L,\Lb)$. A typical picture (if $\ub \geq 1$
) of $D_{\ub,u}$ is the grey region. If $\ub < 1$, the picture of $D_{\ub,u}$ looks like a triangle:

\includegraphics[width = 4.5 in]{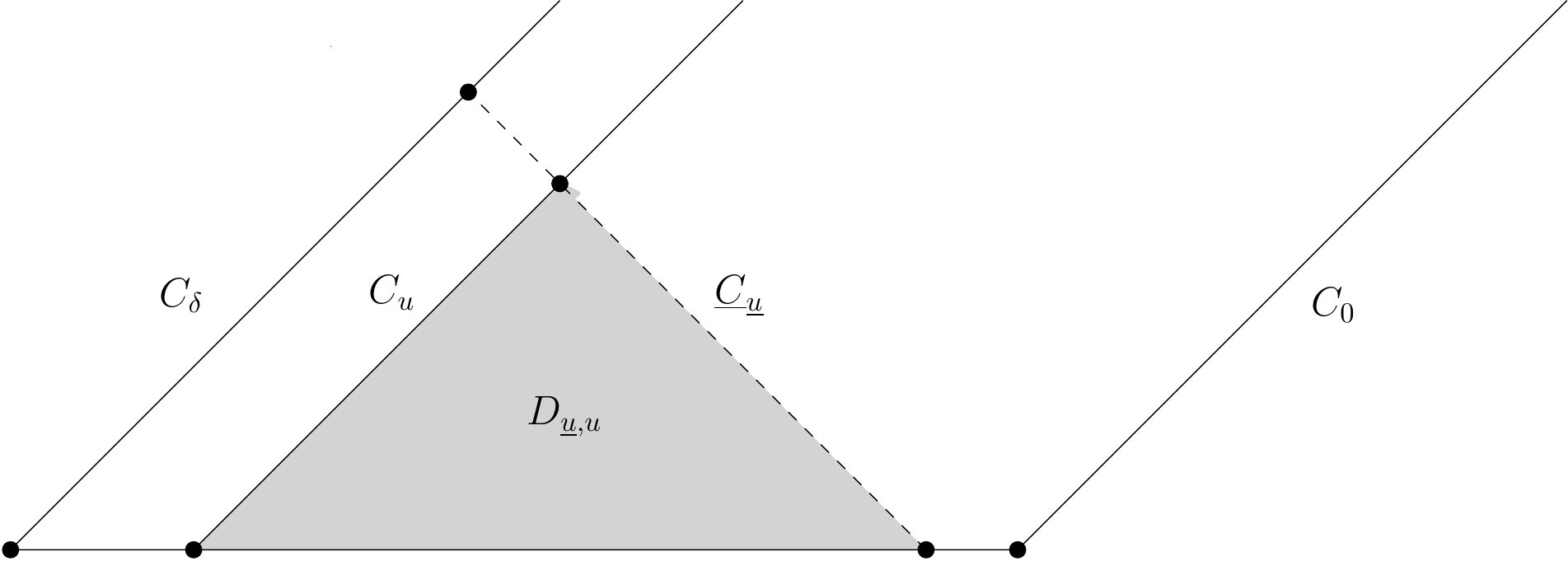}

We remark that, for both cases, both $\{C_{u'} | 0 \leq u' \leq u \}$ and $\{\Cb_{\ub'} | 1-u\leq \ub' \leq \ub\}$ foliate $D_{\ub,u}$.

\bigskip

In view of Remark \ref{geometric interpretation of the short pulse data}, we also remark that the choice of the short pulse data is also adapted to the double null foliation in the short pulse region: the data is chosen in a way that very little energy propagates in the incoming direction through $C_{u}$'s. We expect most of the energy will radiate through the $\Cb_{\ub}$'s to the future null infinity.
\bigskip

For a given $2$-sphere $S_{\ub,u}$, we use $\gslash$ to denote the induced metric from the Minkowski metric on $S_{\ub,u}$. The intrinsic covariant derivative on $S_{\ub,u}$ is denoted by $\nablaslash$. This covariant derivative is closely related to the rotational symmetry of $\mathbb{R}^{3+1}$. Recall that the infinitesimal rotations are represented by the following three vector fields:
\begin{equation*}
 \Omega_{ij}=x_i\partial_j -x_j\partial_i, \ \ \text{for }  1\leq i<j\leq 3.
\end{equation*}
We use  $\Omega$ as a short hand notation for an arbitrary choice from the above vector fields. We also use $\Omega^2$ to denote an operator of the form $\Omega_{i'j'}\Omega_{ij}$; similarly for $\Omega^n$. For a given function $\phi$, we use $|\Omega \phi|$ to denote $\sum |\Omega_{ij} \phi|$ and use $|\Omega^2 \phi|$ to denote $\sum |\Omega_{ij} \Omega_{i'j'}\phi|$, and so on. Therefore, by a direct computation, we obtain
\begin{equation*}
|\Omega \phi| \sim r|\nablaslash \phi|,
\end{equation*}
where the size $\nablaslash \phi$ is measured with respect to $\gslash$. Moreover, for all $n$, we have
\begin{equation*}
|\Omega^n \phi| \sim_n  r^n|\nablaslash^n \phi|.
\end{equation*}
In the rest of the paper, the number of derivatives that we impose on the solution is a fixed number which does not exceed, say, $30$. Therefore the dependence on $n$ in the above inequality is universal.

\begin{remark}
In the short pulse region II, if $\delta$ is sufficiently small, then $|\ub| \sim r$. Therefore, for all $n$, we have
\begin{equation*}
|\Omega^n \phi| \sim |\ub|^n |\nablaslash^n \phi|.
\end{equation*}
In particular, for all $p\geq 1$  and $n$, we have
\begin{equation*}
\|\Omega^n \phi\|_{L^p(S_{\ub,u})} \sim |\ub|^n \|\nablaslash^n\phi\|_{L^p(S_{\ub,u})}.
\end{equation*}
\end{remark}

\bigskip

\subsection{Comments on the proof}
We construct a solution $\phi$ in three steps:
\begin{itemize}
\item \underline{Step 1} We construct $\phi$ in the 
short pulse region II.

The initial data for this region are given on initial hypersurface 
$\Sigma_1$. 
We expect to see the largeness of the data in the proof. In particular, the $\Lb$ derivative of the solution causes a loss of $\delta^{-1}$. This makes the proof difficult and also different from the classical small data problem. The decay of derivatives of $\phi$ is another difficulty which will be explained in detail in Section \ref{Main features of the proof}.

\item \underline{Step 2} Smallness of $\phi$ on $C_\delta$.

Although $\phi$ constructed in \underline{Step 1} has large derivatives, we show that the derivatives of $\phi$ are indeed small on the inner boundary $C_\delta$.  This is a key intermediate step: since in next step, the $\phi$ restricted on $C_\delta$ gives initial characteristic data, this step allows one to reduce the problem to a small data problem in region I.

\item \underline{Step 3} We construct $\phi$ in the small data region 
I.

In region I, the problem is reduced to a small data problem. We can then use the classical approach to construct $\phi$.
\end{itemize}

\subsubsection{Vector field method}
We will derive energy estimates for the main equation \eqref{Main Equation}. Our approach is based on the classical vector field method and we briefly recall the main structure of the method as follows.

Let $\phi$ be a (scalar) solution for a non-homogenous wave equation $\Box \phi = \Phi$ on $\mathbb{R}^{3+1}$. The energy-stress tensor associated to $\phi$  is $ \mathbb{T}_{\alpha \beta}[\phi] = \nabla_\alpha \phi \nabla_\beta \phi -\frac{1}{2}g_{\alpha \beta} \nabla^{\mu} \phi \nabla_{\mu} \phi$ where $g_{\alpha\beta}$ is the Minkowski metric. In particular, in terms of the null pair $(L,\Lb)$, we have $\mathbb{T}[\phi](L,L)=(L\phi)^{2}$, $\mathbb{T}[\phi](\underline{L},\underline{L})=(\underline{L}\phi)^{2}$ and $\mathbb{T}[\phi](L,\underline{L})
=|\slashed{\nabla}\phi|^{2}$. Given a vector field $X$, we use $ ^{(X)}\pi_{\mu \nu} = \frac{1}{2} \mathcal{L}_X g_{\mu \nu}$ to denote its deformation tensor. The energy currents associated to $\phi$ are defined by $J^{X}_{\alpha}[\phi] = \mathbb{T_{\alpha\mu}}[\phi]X^{\mu}$ and $K^X [\phi] = \mathbb{T}^{\mu\nu}[\phi]\, ^{(X)}\pi_{\mu \nu}$. The following divergence identity is the key to the energy estimates:
\begin{equation}\label{divergence of J}
 \nabla^\alpha J^{X}_{\alpha}[\phi] =  K^X [\phi] + \Phi \cdot X \phi.
\end{equation}

In applications, we integrate this identity on the spacetime region. This is equivalent to multiplying $\Box \phi = \Phi$ by $X\phi$ and then integrating by parts. This is the reason that we call $X$ a \emph{multiplier} vector field.

\bigskip

In the short pulse region II, we integrate \eqref{divergence of J} on $D_{\ub,u}$. Since $\phi$ and its derivatives vanish on $C_0$, this yields
\begin{equation}\label{fundamental energy identity}
 \int_{C_u} \mathbb{T}[\phi](X,L)+\int_{\Cb_{\ub}} \mathbb{T}[\phi](X,\Lb)  = \frac{1}{2} \int_{\Sigma_{1}} \mathbb{T}[\phi](X,L+\Lb) +\doubleint_{D_{\ub,u}}  K^X [\phi] + \Phi \cdot X \phi.
\end{equation}
where $\Sigma_1$ is the initial Cauchy hypersurface $\{t=1\}$.

\bigskip

In the short pulse region II, we will use two multiplier vector fields: $X = \Lb$ and $X = \ub^{\alpha} L$, where the power $\alpha =1 -\varepsilon_0$ and $\varepsilon_0 \in (0,\frac{1}{2})$ is a given constant. The first plays a similar role to the time vector field $\partial_t$; the second plays a similar role to the vector field $S=t\partial_t+r\partial_r$ as a multiplier vectorfield.

For  $X = \Lb$ and $X = \ub^{\alpha} L$, the corresponding deformation tensors and energy currents are
\begin{itemize}
\item For $X=\Lb$, $ \pi_{AB} = -\frac{1}{r}\gslash_{AB}$ and $ K =-\frac{1}{r}L\phi\cdot\Lb \phi$.

\item For $X=\ub^{\alpha} L$, $\pi_{L \Lb}=-\alpha \ub^{\alpha-1}$, $\pi_{AB} = \frac{1}{r}\ub^{\alpha} \gslash_{AB}$ and $K=-\frac{\alpha}{2}\ub^{\alpha-1}|\nablaslash\phi|^{2}+\frac{1}{r}\ub^{\alpha}L\phi\cdot\Lb \phi$.
\end{itemize}
respectively. We remark that indices $A$ and $B$ are used to denote a frame on $S_{\ub, u}$ and we only listed the nonzero components of the deformation tensors.

\bigskip

We will also need estimates for higher order derivatives for $\phi$. To achieve this, we will commute the main equation \eqref{Main Equation} with certain vector fields, i.e. the \emph{commutator vector fields}. These vector fields are essentially the Lie algebras of the conformal isometries of $\mathbb{R}^{3+1}$. We list all of them as follows:
\begin{equation*}
\mathcal{Z} = \big\{\Omega_{ij}, \Omega_{0i}, \partial_t, \partial_i, S \big| i,j = 1,2,3, i\neq j\big\},
\end{equation*}
where $\Omega_{0i}=x_{i}\partial_{t}+t\partial_{i}$ and $S=t\partial_{t}+r\partial_{r}=\ub L+u\Lb$. We also define the \emph{good} and \emph{bad} commutator vector fields:
\begin{equation*}
\mathcal{Z} = \mathcal{Z}_g \sqcup \mathcal{Z}_b,\ \ \mathcal{Z}_b = \big\{\partial_t, \partial_i \big|i=1,2,3 \big\}.
\end{equation*}
As shorthand notations, we use $Z$ to denote an arbitrary vector field from $\mathcal{Z}$; similarly, we use $Z_g$ and $Z_b$ to denote vectors from $\mathcal{Z}_g$ and $\mathcal{Z}_b$ respectively. Geometrically, a good vector field $Z_g$ is tangential to the outgoing light cone $C_0$, but a bad vector field $Z_b$ is transversal to $C_0$.

\subsubsection{A word on null forms}
Recall that a quadratic form $Q$ over $\mathbb{R}^{3+1}$ is a \emph{null form} if $Q(\xi, \xi) = 0$ for all null vector $\xi \in \mathbb{R}^{3+1}$. The space of null forms are spanned by the following seven forms: $Q_0 (\xi, \eta) = g(\xi, \eta)$ and $Q_{\alpha \beta}(\xi, \eta) = \xi_\alpha \eta_\beta - \eta_\alpha \xi_\beta (0\leq \alpha, \beta \leq 3)$. Given scalar functions $\phi$, $\psi$ and a null form $Q(\xi,\eta) = Q^{\alpha\beta}\, \xi_\alpha \, \eta_\beta$, we use $Q(\nabla \phi, \nabla \psi)$ as a shorthand for $Q(\nabla \phi, \nabla \psi) = Q^{\alpha\beta} \, \partial_\alpha \phi \, \partial_\beta \psi.$

\bigskip

For a (conformal) Killing vector field $Z \in \mathcal{Z}$, we have
\begin{equation}\label{null form commuatator}
 Z Q(\nabla \phi, \nabla \psi) = Q(\nabla Z \phi, \nabla \psi) + Q(\nabla \phi, \nabla Z \psi) + \widetilde{Q}(\nabla \phi, \nabla \psi),
\end{equation}
where $\widetilde{Q}$ is a null form, which may or may not be $Q$.

\bigskip

In terms of the null pair $(L,\Lb)$, a null form $Q$ satisfies the following pointwise estimates
\begin{equation}\label{null form bound}
 |Q(\nabla \phi, \nabla \psi)| \lesssim |L \phi| \,|\Lb\psi| + |\Lb \phi|\,|L\psi| + |\nablaslash \phi|\, |\nablaslash \psi| +\big(| L\phi| + |\Lb \phi|\big)|\nablaslash \psi|  + |\nablaslash \phi|\big(|L\psi| +  |\Lb\psi|\big) .
\end{equation}
In particular, on the right hand side of the inequality, the term $|\Lb \phi|^2$ does not appear.

\subsubsection{Main features of the proof}\label{Main features of the proof}

We discuss main difficulties of the problem and also the ideas to get around them.

\bigskip

\begin{itemize}
\item Largeness/Loss of $\delta^{-1}$ in the short pulse region.

\bigskip

In the short pulse region, if one differentiates $\phi$ in the $\Lb$ direction, then the resulting function will be approximately $\delta^{-1}$ times as large as the initial functions. Schematically, we can regard $\Lb$ as $\Lb \sim \delta^{-1}$. Similarly, $L \sim 1$ and $\nablaslash \sim 1$.

The large factor $\delta^{-1}$ maybe fatal to the energy estimates for nonlinear terms. The resolution of this difficulty is exactly the basic philosophy of null conditions: if
one term behaves badly, say 
$|\Lb \phi| \sim \delta^{-\frac{1}{2}}$ in the nonlinearities,
it must be coupled with the a good term, say $L\phi$ or $\nablaslash \phi$, which are both of size $\delta^{\frac{1}{2}}$.
Their product will then be a term of size $1$ which will be manageable in the proof. \bigskip

\item Relaxation in $\delta$ for the propagation estimates.

\bigskip

On the initial hypersurface $\Sigma_{1}$, it is easy to see that the data satisfy $\|\Lb \phi\|_{L^{\infty}_{\Sigma_{1}}} \sim \delta^{-\frac{1}{2}}$ and $\|\nablaslash \phi\|_{L^{\infty}_{\Sigma_{1}}}\sim \delta^{\frac{1}{2}}$. Up to a correct decay factor in $t$, we hope the size of $\Lb \phi$ and $\nablaslash \phi$ measured in $\delta$ can be propagated for later times, i.e. $\|\Lb \phi\|_{L^{\infty}_{\Sigma_{1}}} \sim \delta^{-\frac{1}{2}}$ and $\|\nablaslash \phi\|_{L^{\infty}_{\Sigma_{1}}}\sim \delta^{\frac{1}{2}}$ should be always true. Recall that the proof will be based on energy estimates. If we use $\Lb$ as a multiplier vector field and integrate in $D_{\ub, u}$ in the short pulse region, the energy on the left hand side of \eqref{fundamental energy identity} is $\int_{\Cb_{\ub}}|\Lb \phi|^{2}+ \int_{C_{u}}|\nablaslash \phi|^{2}$. Therefore, the expected propagation estimates suggest that $\int_{\Cb_{\ub}}|\Lb \phi|^{2} \lesssim 1$ and $\int_{C_{u}}|\nablaslash \phi|^{2} \lesssim \delta$. Therefore, in view of the form of the
energy, the disparity of the $\delta$ power for these two quantities only gives the desired bound for $\Lb\phi$, but not for $\nablaslash \phi$. This may lead to the failure of closing the bootstrap argument.

To get around this difficulty, we pretend that the amplitude of $\nablaslash \phi$ was worse than that suggested by the initial data. The purpose of this relaxation is to make the two terms in $\int_{\Cb_{\ub}}|\Lb \phi|^{2}+\int_{C_{u}}|\nablaslash \phi|^{2}$ comparable. More specifically, we first prove that $\|\nablaslash \phi\|_{L^{\infty}(\Sigma_{1})}\lesssim 1$ can be propagated by controlling $\int_{\Cb_{\ub}}|\Lb \phi|^{2}+\int_{C_{u}}|\nablaslash \phi|^{2}$ in terms of initial energies. In this way we are able to close the bootstrap argument. Then we use the proved energy estimate to recover the estimate $\|\nablaslash\phi\|_{L^{\infty}(\Sigma_{t})}\lesssim\delta^{1/4}t^{-2}$ by affording to lose one derivative.


\bigskip

\item Relaxation in the decay factor in the short pulse region.

\bigskip

According to the decay rate of linear waves, one may expect the decay of $\phi$ or more precisely the derivatives of $\phi$ should be $\frac{1}{t}$ or $\frac{1}{\ub}$ in the short pulse region. This expected decay will cause a loss of $\log t$ in the energy estimates since we may need to integrate a factor of size $\frac{1}{t}$ coming from the nonlinear term.

The idea to get around this point is also to relax the decay rate a little bit. This is why we choose $X = |\ub|^{1-\varepsilon_0}L$ as a multiplier vector field instead of using the standard vectorfield $S$. The $|\ub|^{-\varepsilon_0}$ will be amplified to $|\ub|^{-2\varepsilon_0}$ in the energy estimates due to the nonlinearity. Therefore, we can gain a little more decay relative to the relaxed decay. This is just enough to close the argument for the energy estimates.

\bigskip

\item Smallness of the solution on $C_\delta$.

\bigskip

This is precisely the question that we will answer in Step 2 of the proof. As we discussed, in the short pulse region, we expect $\Lb \sim \delta^{-1}$. In particular, 
we expect that, for all the bad vector fields $Z_b$, we also have $Z_b \sim \delta^{-1}$. Therefore, for a given $n$, the restriction of $Z_b{}^n \phi$ on $C_\delta$ may be of size $\delta^{\frac{1}{2}-n}$. This is by no means small.

The key point of the proof is the following observation: on the 2-sphere $S_{1-\delta,\delta}$, i.e. the initial sphere of $C_\delta$, the data vanish completely since they are compactly supported on $\Sigma_1$ between $S_{1-\delta,\delta}$ and $S_{1,0}$. Therefore, even the bad derivatives of $\phi$ are small initially. To get the smallness of $\phi$, we will integrate along null geodesics on $C_{\delta}$ to trace all the information back to the data. In this way, we can show that up to an error of size $\delta^\frac{1}{2}$, all derivatives of $\phi$ are comparable to their initial values.

\end{itemize}

\subsubsection{Comparison with the previous work \cite{W-Y-12}}
We now discuss the difference between the present work and the previous work \cite{W-Y-12}.\\

\begin{itemize}
\item Cauchy data versus characteristic scattering data.\\

In the present work, we consider the Cauchy data given on $\Sigma_{1}$ satisfying \eqref{data I} and \eqref{data II}. However \eqref{data I}, \eqref{data II} do not give the a priori decay profile. To prove pointwise decay estimates for $\phi$, besides the standard vectorfield $\Lb$, one has to choose $\ub^{\alpha}L, \alpha\in(1/2,1)$ as a multiplier vectorfield. Here the index $\alpha$ is chosen in such a way that it is enough to prove decay estimates but the decay rate is not too strong to prove.

This should be compared with the characteristic scattering data considered in \cite{W-Y-12} given at the past null infinity $C_{-\infty}$. More specifically, the data for $\phi$ in \cite{W-Y-12} has the following form:
\begin{align}\label{characteristic data}
\lim_{u\rightarrow-\infty}|u|\phi(u,\ub,\theta)=\delta^{1/2}\psi_{0}\left(\frac{\ub}{\delta},\theta\right)
\end{align}
Here $\psi_{0}:(0,1)\times \mathbb{S}^{2}\rightarrow\mathbb{R}$ is a compactly supported smooth function. Note that the data \eqref{characteristic data} has the property

\begin{align}\label{scatter}
\phi(u,\ub,\theta)\sim\frac{\delta^{1/2}}{|u|}\psi_{0}\left(\frac{\ub}{\delta},\theta\right)+o\left(\frac{1}{|u|}\right), \quad u\rightarrow-\infty.
\end{align}
Therefore the decay profile $\frac{1}{|u|}$ for $\phi$ is given a priori, which simplifies the proof of pointwise decay estimates for $\phi$ and its derivatives. In fact, the authors in \cite{W-Y-12} use the standard vectorfields $L$ and $\Lb$ as the multipliers to prove the energy estimates.\\

\item Compactly supported data versus non-compactly supported trace.\\

In \cite{W-Y-12} the characteristic data \eqref{characteristic data} given at the past null infinity is compactly supported in $u\in(0,\delta)$. After solving the characteristic problem all the way up to $t=-1$, the restriction of solution on $\Sigma_{-1}$ gives the data for the time-reversed Cauchy problem. The Cauchy data given in this way is implicitly and can never be compactly supported. On the contrary, the Cauchy data in the present work is given directly and explicitly. Moreover, as it is shown in \eqref{data example}, the data can be compactly supported. On the other hand, compared to the characteristic data at the past null infinity, which is compactly supported, the trace of the solution at the future null infinity ($\ub=\infty$) in the present work is not compactly supported.\\

\end{itemize}

\subsubsection{Applications in physical problems}
We would like to discuss further applications of our method to other wave type equations, especially Yang-Mills equations in gauge theory and Einstein equations in general relativity. For both systems of equations, there is no known result to derive global asymptotic behaviors for large data problem. Taking Yang-Mills equations as an example. Let $F$ be the Yang-Mills field. We define two 1-forms by contracting with $L$ and $\Lb$ on $S_{u,\ub}$: $\alpha_F = i_L F$, $\underline{\alpha}_F = i_{\Lb} F$. Since $F$ is an Lie-algebra valued two forms (of dimension $6$), the rest two components of $F$ are denoted by $\rho_F$ and $\sigma_F$. This four components $\alpha_F$,$\underline{\alpha}_F$, $\rho_F$ and $\sigma_F$ consists of a complete decomposition of $F$ by using null frames. To make connections to our method, we make the following correspondence: 
$$\alpha_F \mapsto L\phi, \ \ \underline{\alpha}_F \mapsto \Lb \phi, \ \ (\rho_F,\sigma_F)\mapsto \nablaslash \phi.$$
By the correspondence, we can use the corresponding energy ansatz for the each component respectively. Since Yang-Mills equations also have null structures, we expect our method can prove the first asymptotic description large data problem. Similarly, we can also study Einstein equations in such a way (by virtue of harmonic coordinates). This will be a forthcoming paper.

\subsubsection{Outline of the paper}



The rest of the paper is organized as follows:

In Section \ref{Short pulse region}, we establish \emph{a priori} energy estimates for higher order derivatives of the solution in the short pulse region. As consequences, first of all, we can construct the solution in the short pulse region; Secondly, we can obtain a smallness estimate for the solution on $C_{\delta}$, i.e. the inner boundary of the short pulse data region.

In Section 3, with a modified Klainerman-Sobolev inequality, we construct global solutions in the small data region.

\section{Short pulse region}\label{Short pulse region}

The goal of the current section is to construct the solution $\phi$ in the short pulse region. The construction relies on a priori energy estimates. We assume that the solution $\phi$ exists on spacetime domain $D_{\ub^*,u^*}$. This domain is inside the short pulse region, i.e. $u^*\in (0,\delta)$ and $\ub^* \in (1-u^*,+\infty)$.

\bigskip

We first introduce the energy norms. Let $u,u' \in (0,u^*)$ and $\ub,\ub' \in (1-u^*, \ub^*)$. Let $C^{\ub'}_u$ be the part of the cone $C_u$ so that $1-u^* \leq \ub \leq \ub'$ and let ${\Cb}^{u'}_{\ub}$ be the part of the cone ${\Cb}_{\ub}$ so that $0 \leq u \leq u'$. Whenever there is no confusion, we will use $C_u$ and $\Cb_\ub$ instead of  $C^{\ub}_u$ and $\Cb^u_{\ub}$. We use $\Sigma_{1}$ to denote the annulus region 
$\{(r,\theta)\mid 1-2\delta \leq r \leq 1\}$ on $\{t=1\}$ in the current section.

For a given $k \in \mathbb{Z}_{\geq 0}$, we introduce the following homogeneous norms:
\begin{equation*}
\begin{split}
E_k(u,\ub) &= \sum_{Z_g \in \mathcal{Z}_g, Z_b\in\mathcal{Z}_b \atop 0\leq l\leq k} \Big( \delta^{l}\|\slashed{\nabla}Z_{b}^{l}Z_{g}^{k-l}\phi\|_{L^2(C^{\ub}_u)}+\delta^{l-\frac{1}{2}}\||\ub|^{\frac{\alpha}{2}}LZ_{b}^{l}Z_{g}^{k-l}\phi\|_{L^2(C^{\ub}_u)} \Big),\\
\Eb_k(u, \ub) &= \sum_{Z_g \in \mathcal{Z}_g, Z_b\in\mathcal{Z}_b \atop 0\leq l\leq k} \Big(\delta^{l}\|\Lb Z_{b}^{l}Z_g^{k-l} \phi\|_{L^2(\Cb^u_{\ub})}  +\delta^{l-\frac{1}{2}}\||\ub|^{\frac{\alpha}{2}}\slashed{\nabla}Z_{b}^{l}Z_{g}^{k-l}\phi\|_{L^2(\Cb^u_{\ub})}\Big).
\end{split}
\end{equation*}
We also introduce the inhomogeneous norms:
\begin{equation*}
E_{\leq k} (u,\ub) = \sum_{0\leq j \leq  k} E_{j} (u,\ub), \ \ \Eb_{\leq k}(u, \ub) = \sum_{0\leq j \leq  k} \Eb_{j} (u,\ub).
\end{equation*}

On the initial hypersurface $\Sigma_1$, we introduce the following initial energy norms:
\begin{equation*}
E_{\leq n}( \Sigma_1)= \sum_{Z_g \in \mathcal{Z}_g, Z_b\in\mathcal{Z}_b \atop 0\leq l\leq k, 0\leq k \leq n} \delta^{l}\|\Lb Z_{b}^{l}Z_{g}^{k-l}\phi\|_{L^{2}{(\Sigma_{1})}} + \delta^{l-\frac{1}{2}} \|\nablaslash Z_{b}^{l}Z_{g}^{k-l} \phi\|_{L^{2}{(\Sigma_{1})}}+\delta^{l-1}\|L Z_{b}^{l}Z_{g}^{k-l}\phi\|_{L^{2}{(\Sigma_{1})}}.
\end{equation*}
According to the behavior of $\phi_{0}$ and $\phi_{1}$ on $\Sigma_{1}$ as well as the properties of commutation vectorfields, using the arguments in \cite{W-Y-12} and \cite{W-Y-13}, we have the following lemma: 
\begin{lemma}
For all $\delta$, we have
\begin{equation}\label{initial bound}
E_{\leq n-1}( \Sigma_1) \lesssim I_n.
\end{equation}
Here $I_{n}$ is a constant depending only on $n$.
\end{lemma}

\subsection{Main a priori estimates}
This subsection is the central part of the paper. The goal is to bound $E_{\leq 3}(u,\ub)$ on $D_{\ub^*,u^*}$ where the solution $\phi$ is assumed to exist.
\begin{proposition}\label{proposition bootstrap}
There exists $\delta_0 >0$, so that for all $\delta < \delta_0$,  for all $u \in (0,u^*)$ and $\ub \in (1-u^*, \ub^*)$, we have
\begin{equation}\label{main estimates}
E_{\leq 3}(u,\ub) + \Eb_{\leq 3}( u,\ub) \leq C(I_4),
\end{equation}
where $C(I_4)$ is a constant depending only on $I_4$.
\end{proposition}

The proof of the proposition is based on a standard bootstrap argument. On $D_{\ub^*,u^*}$, since we assume that $\phi$ exists, there is a large constant $M$, so that
\begin{equation}\label{bootstrap assumption}
E_{\leq 3}(u,\ub) + \Eb_{\leq 3}( u,\ub) \lesssim M,
\end{equation}
for all $u \in (0,u^*)$ and $\ub \in (1-u^*, \ub^*)$. The large constant $M$ may depend on $\phi$ itself. The purpose of the bootstrap argument is to show that, if $\delta$ is sufficiently small, then one can choose $M$ in such a way that it depends only on $I_4$. Hence, we obtain the proof of \eqref{main estimates}

\subsubsection{Preliminary estimates}
The goal of this subsection is to use the bootstrap assumption \eqref{bootstrap assumption} to get estimates on lower order derivatives of $\phi$ (up to second derivatives).

We first recall the Sobolev inequalities on $S_{\ub,u}$, $\Cb_{\ub}$ and $C_u$ in the short pulse region. Recall that in the short pulse region, we have $|\ub|\sim r$ provided $\delta$ is sufficiently small. Let $\phi$ be a smooth function.

On $S_{\ub,u}$, we have
\begin{equation}\label{L-infty Sobolev}
\begin{split}
\|\phi\|_{L^\infty(S_{\ub,u})}&\lesssim |\ub|^{-1/2}\big(\|\phi\|_{L^{4}(S_{\ub,u})}+\|\Omega\phi\|_{L^{4}(S_{\ub,u})}\big),
\end{split}
\end{equation}
\begin{equation}\label{L4 Sobolev}
\|\phi\|_{L^{4}(S_{\ub,u})}\lesssim |\ub|^{-1/2}\big(\|\phi\|_{L^{2}(S_{\ub,u})}+\|\Omega\phi\|_{L^{2}(S_{\ub,u})}\big).
\end{equation}

On $\Cb_\ub$, if in addition we assume that $\phi \equiv 0$ on $C_0$, we have
\begin{equation}\label{L2 Sobolev}
\begin{split}
\|\phi\|_{L^{2}(S_{\underline{u},u})}&\lesssim\|\underline{L}\phi\|^{1/2}_{L^{2}(\underline{C}_{\underline{u}})}\|\phi\|^{1/2}_{L^{2}(\underline{C}_{\underline{u}})},\\
\|\phi\|_{L^{4}(S_{\underline{u},u})}&\lesssim |\ub|^{-\frac{1}{2}}\|\underline{L}\phi\|^{1/2}_{L^{2}(\Cb_{\ub})}\big(\|\phi\|^{1/2}_{L^{2}(\underline{C}_{\underline{u}})}+\|\Omega\phi\|^{1/2}_{L^{2}(\underline{C}_{\underline{u}})}\big).
\end{split}
\end{equation}
We remark that the assumption $\phi \equiv 0$ on $C_0$ will be always true when we apply the above inequalities in the short pulse region in the rest of the paper, since the solution $\phi$ of the main equations \eqref{Main Equation} (if it exists) vanishes to infinite order on $C_0$.

If $\phi$ is supported in the annular region $\{(r,\theta)|1-\delta\leq r\leq 1\}$ on the initial Cauchy hypersurface $\Sigma_1$, we have
\begin{equation}\label{L2 Sobolev initial}
\|\phi\|_{L^{2}(S_{1-u,u})}\lesssim \delta^{1/2}\big(\|\Lb\phi\|_{L^{2}(\Sigma_{1})}+\|L\phi\|_{L^{2}(\Sigma_{1})}\big).
\end{equation}

For the proof of the above inequalities, we refer the reader to \cite{Ch-08}.

\bigskip

We also recall the Gronwall's inequality. Let $f(t)$ be a non-negative function defined on an interval $I$ with initial point $t_{0}$. If $f$ satisfies
\begin{equation*}
 \frac{d}{dt}f\leq a\cdot f+b
\end{equation*}
with two non-negative functions $a, b\in L^{1}(I)$, then for all $t\in I$, we have
\begin{equation*}
 f(t)\leq e^{A(t)}\big(f(t_{0})+\int_{t_{0}}^{t}e^{-A(\tau)}b(\tau)d\tau\big)
\end{equation*}
where $A(t)=\int_{t_{0}}^{t}a(\tau)d\tau$.

\bigskip

We start to derive estimates and we treat $\ub$ as a fixed constant. By virtue of null pair $(L,\Lb)$, we rewrite the main system of equations \eqref{Main Equation} as
\begin{equation}\label{Equation in null frame}
 -L\Lb\phi+\slashed{\triangle}\phi+\frac{1}{r}(L\phi-\Lb\phi)=Q(\nabla \phi,\nabla \phi).
\end{equation}
We remark that, $\phi$ is now a $\mathbb{R}^N$-valued function and the norms used in the rest of the paper are with respect to a fixed inner product in $\mathbb{R}^N$. For example, the symbol $|L\phi|$ denotes $\sqrt{\sum_{I\leq N}(L\phi^{I})^{2}}$.

We also need to commute derivatives with \eqref{Equation in null frame}. Recall that, for all $Z \in \mathcal{Z}$ except for $Z = S$, we have $[\Box, Z]=0$. Indeed, we have $[\Box, S] = 2\Box$. Combining this remark with \eqref{null form commuatator}, for all $k \geq 0$, we can commute $k$ vectors $Z_1,Z_2,\cdots, Z_k \in \mathcal{Z}$ with \eqref{Main Equation} to obtain a semilinear wave equation for $Z_1 Z_2 \cdots Z_k \phi$. We use the shorthand notation $Z^k\phi$ to denote $Z_1 Z_2 \cdots Z_k \phi$, therefore, we have
\begin{equation}\label{commuted null frame equation}
\Box Z^{k} \phi = \sum_{p+q \leq k}Q(\nabla Z^p \phi, \nabla Z^q \phi).
\end{equation}

\bigskip

We combine \eqref{L-infty Sobolev}, \eqref{L4 Sobolev}, \eqref{L2 Sobolev} and bootstrap assumption \eqref{bootstrap assumption}. We first have
\begin{align*}
 \| \nablaslash \phi \|_{L^{4}(S_{\ub,u})} &\lesssim |\ub|^{-\frac{1}{2}} \| \Lb \nablaslash\phi\|^{\frac{1}{2}}_{L^{2}(\Cb_{\ub})}\big(\| \nablaslash \phi \|^{\frac{1}{2}}_{L^{2}(\Cb_{\ub})} +|\ub|^{\frac{1}{2}}\| \nablaslash^2 \phi \|_{L^{2}(\Cb_{\ub})}\big)\\
 &\lesssim  |\ub|^{-\frac{1}{2}}(|\ub|^{-1}M )^\frac{1}{2}\big((\delta^{\frac{1}{2}}|\ub|^{-\frac{\alpha}{2}}M)^\frac{1}{2}+ |\ub|^\frac{1}{2}(\delta^{\frac{1}{2}} |\ub|^{-\frac{2+\alpha}{2}}M)^\frac{1}{2}\big).
\end{align*}
Hence,
\begin{equation*}\label{p_4}
\| \nablaslash\phi\|_{L^{4}(S_{\ub,u})}\lesssim \delta^{\frac{1}{4}}|\ub|^{-1-\frac{\alpha}{4}}M.
\end{equation*}

Similarly, since in the bootstrap assumption \eqref{bootstrap assumption}, we have assumed bounds on four derivatives on $\phi$, we can repeat the above argument to derive
\begin{equation}\label{p_5}
\| \nablaslash\Omega \phi\|_{L^{4}(S_{\ub,u})}\lesssim \delta^{\frac{1}{4}}|\ub|^{-1-\frac{\alpha}{4}}M,
\end{equation}
and
\begin{equation}\label{p_6}
\| \nablaslash\Omega^2 \phi\|_{L^{4}(S_{\ub,u})}\lesssim \delta^{\frac{1}{4}}|\ub|^{-1-\frac{\alpha}{4}}M.
\end{equation}

Combining \eqref{p_4} and \eqref{p_5}, Sobolev inequality implies
\begin{equation}\label{p_7}
\begin{split}
 \|\nablaslash \phi\|_{L^\infty(S_{\ub,u})} &\lesssim  |\ub|^{-\frac{1}{2}}\big(\| \nablaslash\phi \|_{L^{4}(S_{\ub,u})}
 +\| \nablaslash\Omega \phi \|_{L^{4}(S_{\ub,u})}\big)\\
 &\lesssim \delta^{\frac{1}{4}}|\ub|^{-\frac{3}{2}-\frac{\alpha}{4}}M.
 \end{split}
 \end{equation}

Similarly, we have
\begin{equation}\label{p_8}
 \|\nablaslash\Omega \phi\|_{L^\infty(S_{\ub,u})} \lesssim  \delta^{\frac{1}{4}}|\ub|^{-\frac{3}{2}-\frac{\alpha}{4}}M.
\end{equation}

By repeating the above argument, for $0 \leq l \leq k \leq 2$, we can also easily obtain
\begin{equation*}
\| \nablaslash Z_b^l Z_g^{k-l} \phi\|_{L^{4}(S_{\ub,u})}\lesssim \delta^{\frac{1}{4}-l}|\ub|^{-1-\frac{\alpha}{4}}M.
\end{equation*}
and for $0 \leq l \leq k \leq 1$
\begin{equation}\label{p_9}
 \|\nablaslash Z_b^l Z_g^{k-l}\phi\|_{L^\infty(S_{\ub,u})} \lesssim  \delta^{\frac{1}{4}-l}|\ub|^{-\frac{3}{2}-\frac{\alpha}{4}}M.
\end{equation}

\bigskip

We turn to the bound of $L\phi$ in $L^\infty(S_{\ub,u})$. Let $a=\frac{1}{r}+|\Lb\phi|+|\slashed{\nabla}\phi|$ and $b=|\slashed{\triangle}\phi|+\frac{1}{r}|\Lb\phi|+|\Lb\phi\slashed{\nabla}\phi|+|\slashed{\nabla}\phi|^{2}$, in view of \eqref{null form bound}, \eqref{Equation in null frame} yields
\begin{equation*}
\Lb|L\phi|\lesssim a|L\phi|+b.
\end{equation*}
We would like to integrate this equation directly along $\Lb$ to derive the pointwise bound on $L\phi$. Since $L\phi$ vanishes along $C_{0}$, in view of Gronwall's inequality, it suffices to control $\|a\|_{L^{1}_{u}L^{\infty}(S_{\underline{u},u})}$ and $\|b\|_{L^{1}_{u}L^{\infty}(S_{\underline{u},u})}$. We only give the estimates on $|\Lb\phi|$ appearing in $a$ and $b$. The others can be estimated directly from \eqref{p_9}. According to Sobolev inequality, we have
\begin{align*}
 \|\Lb \phi\|_{L^1_u L^{\infty}(S_{\ub,u})}&\lesssim |\ub|^{-1}\sum_{0\leq j \leq 2} \|\Omega^j \Lb \phi\|_{L^1_u L^{\infty}(S_{\ub,u})}\\
 &\lesssim |\ub|^{-1}\delta^{\frac{1}{2}}\sum_{0\leq j \leq 2} \|\Omega^j \Lb \phi\|_{L^2_u L^{\infty}(S_{\ub,u})}\\
 &\lesssim |\ub|^{-1}\delta^{\frac{1}{2}} M.
\end{align*}
Finally, we can prove
\begin{align*}
 \|a\|_{L^1_u L^{\infty}(S_{\ub,u})}&\lesssim |\ub|^{-1}\delta^{-\frac{1}{2}} M,\\
  \|b\|_{L^1_u L^{\infty}(S_{\ub,u})}&\lesssim |\ub|^{-2}\delta^{-\frac{1}{2}} M.
\end{align*}
Therefore, Gronwall's inequality provides us the following estimates for $L\phi$:
\begin{equation}\label{L-infty for L}
\|L\phi\|_{L^{\infty}(S_{\ub,u})}\lesssim \delta^{1/2}|\ub|^{-2}M.
\end{equation}

\bigskip

By virtue of \eqref{commuted null frame equation} (where $k\leq 2$), we can also bound $L Z_b^l Z_g^{k-l} \phi$ in $L^2(S_{\ub,u})$ in a similar way. Therefore, for $0 \leq l \leq k \leq 2$, we have
\begin{equation}\label{L-2 for L Zk phi}
\|L Z_b^l Z_g^{k-l} \phi\|_{L^{2}(S_{\ub,u})}\lesssim \delta^{1/2-l}|\ub|^{-1}M.
\end{equation}

\bigskip
We turn to the $L^{\infty}(S_{\ub,u})$ estimates on $\Lb \phi$. We start with a computation of $L\big(\ub^{2}(\Lb\phi)^{2}\big)$:
\begin{align*}
L\big(\ub^{2}(\Lb\phi)^{2}\big)&= 2\ub(\Lb\phi)^{2}+2\ub^{2}(\Lb\phi)(L\Lb\phi)\\
&=2\ub(\Lb\phi)^{2}+2\ub^{2}(\Lb\phi)\big[\slashed{\triangle}\phi+\frac{1}{r}(L\phi-\Lb\phi) + Q(\nabla \phi, \nabla \phi)\big]\\
&\lesssim|\ub|^{2}|\Lb\phi|^{2}\big(\frac{2}{r}-\frac{2}{\ub}+|L\phi|+|\slashed{\nabla}\phi|\big)+|\ub|^{2}|\Lb\phi|\big(|\slashed{\triangle}\phi|+\frac{1}{r}|L\phi|+|L\phi||\slashed{\nabla}\phi|+|\slashed{\nabla}\phi|^{2}\big).
\end{align*}
We make the following important observation: in the short pulse region, $|\frac{2}{r}-\frac{2}{\ub}|\lesssim \frac{\delta}{|\ub|^2}$.  Therefore, if we define $y = |\ub||\Lb\phi|$, according to the estimate obtained so far, the previous computation yields
\begin{equation*}
L y^{2}\lesssim \big(\frac{\delta}{|\ub|^{2}}+ \frac{\delta^\frac{1}{4}}{|\ub|^{2}}M\big) y^{2}+ \frac{\delta^\frac{1}{4}}{|\ub|^{2}} M y.
\end{equation*}
We divide both sides of the equation by $y$, thus, we have
\begin{equation*}
L\big(|\ub||\Lb\phi|\big)\lesssim \big(\frac{\delta}{|\ub|^{2}}+ \frac{\delta^\frac{1}{4}}{|\ub|^{2}}M\big) \big(|\ub||\Lb\phi|\big) + \frac{\delta^\frac{1}{4}}{|\ub|^{2}} M.
\end{equation*}
By integrating directly this equation, if $\delta$ is sufficiently small, we obtain
\begin{equation}\label{a1}
\big||\ub||\Lb\phi|(\ub,u,\theta)-C|1-u||\Lb\phi|(1-u,u,\theta)\big| \lesssim  \delta^\frac{1}{4}M.
\end{equation}
where the absolute constant $C$ comes from the use of Gronwall's inequality.
Therefore, according to \eqref{L2 Sobolev initial} and Lemma \ref{initial bound}, we finally obtain
\begin{equation}\label{L infty on Lb phi}
\|\Lb\phi\|_{L^\infty(S_{\ub,u})} \lesssim \frac{\delta^{-\frac{1}{2}}}{|\ub|} I_3 + \delta^\frac{1}{4}|\ub|^{-1}M.
\end{equation}

\bigskip

We remark that the derivation of \eqref{L-infty for L} and \eqref{L infty on Lb phi} depends on not only on the bootstrap assumption \eqref{bootstrap assumption} but also the main equation \eqref{Main Equation}.
We summarize the estimates derived so far as follows:
\begin{equation}\label{preliminary estimates}
\begin{split}
\|\nablaslash Z_b^l Z_g^{k-l} \phi\|_{L^{4}(S_{\ub,u})}&\lesssim \delta^{\frac{1}{4}-l}|\ub|^{-1-\frac{\alpha}{4}}M, \ \ 0 \leq l \leq k \leq 2,\\
 \|\nablaslash \phi\|_{L^\infty(S_{\ub,u})} &\lesssim  \delta^{\frac{1}{4}}|\ub|^{-\frac{3}{2}-\frac{\alpha}{4}}M, \\
 \|L  \phi\|_{L^{\infty}(S_{\ub,u})}&\lesssim \delta^{\frac{1}{2}}|\ub|^{-2}M,  \\
 \|L Z_b^l Z_g^{k-l} \phi\|_{L^{2}(S_{\ub,u})} &\lesssim \delta^{1/2-l}|\ub|^{-1}M, \ \ 0 \leq l \leq k \leq 2,\\
\|\Lb \phi\|_{L^\infty(S_{\ub,u})} &\lesssim \frac{\delta^{-\frac{1}{2}}}{|\ub|} I_3 + \delta^{\frac{1}{4}}|\ub|^{-1}M,
\end{split}
\end{equation}

\begin{remark}
The bootstrap assumptions \eqref{bootstrap assumption} involve relaxed estimates for $\nablaslash \phi$. Roughly speaking, in \eqref{bootstrap assumption} we expect the behavior of $\nablaslash \phi$ with respect to $\delta$ and $\ub$ is approximately $|\ub|^{-1-\frac{\alpha}{2}}M$, i.e.
\begin{equation*}
 \|\nablaslash \phi\|_{L^\infty(S_{\ub,u})} \sim |\ub|^{-1-\frac{\alpha}{2}}\cdot \delta^0 \cdot M.
\end{equation*}
However, the estimates on $\nablaslash \phi$ in \eqref{preliminary estimates} shows that, by affording two more derivatives (via Sobolev inequalities), we can improve the bound on  $\nablaslash \phi$: we get an extra $\delta^\frac{1}{4}$ factor and an extra $\ub^{-\frac{1}{2}+\frac{\alpha}{4}}$ decay factor, i.e.,
\begin{equation*}
 \|\nablaslash \phi\|_{L^\infty(S_{\ub,u})} \sim |\ub|^{-\frac{3}{2}-\frac{\alpha}{4}}\cdot \delta^\frac{1}{4}\cdot M.
\end{equation*}
\end{remark}

\subsubsection{Estimates on $E_{\leq 2}$ and $\Eb_{\leq 2}$}\label{lower order energy estimates}
Recall that for $Z \in \mathcal{Z}$ and $k\geq 0$, we have
\begin{equation}\label{commuated main equatin for energy estimates}
\Box Z^{k} \phi = \sum_{p+q \leq k}Q(\nabla Z^p \phi, \nabla Z^q \phi).
\end{equation}

In Section \eqref{lower order energy estimates}, we fix $k\leq 2$. Let $l\leq k$ be the number of $Z_b$'s appearing in $Z^k$, i.e. $Z^k = Z_b^lZ_g^{k-l}$. We use the vector field method outlined in the introduction to estimate $E_{\leq 2}$ and $\Eb_{\leq 2}$.

\bigskip

In the fundamental energy identity \eqref{fundamental energy identity} and \eqref{commuated main equatin for energy estimates}, we replace $\phi$ by $Z^k \phi$ and take $X = \Lb$ to obtain
\begin{equation*}
\begin{split}
\int_{C_{u}}|\nablaslash Z^k \phi|^{2}&+\int_{\underline{C}_{\underline{u}}}|\Lb Z^k \phi|^{2} =\int_{\Sigma_1}|\nablaslash Z^k \phi|^{2}+|\Lb Z^k \phi|^{2}  + \doubleint_{D_{\ub,u}}Q(\nabla Z^k \phi, \nabla \phi)\Lb Z^k \phi  \\
&  + \sum_{p+q \leq k,\atop p<k, q<k} \doubleint_{D_{\ub,u}}
Q(\nabla Z^{p}\phi, \nabla Z^{q}\phi)\Lb Z^i \phi -\doubleint_{D_{\ub,u}} \frac{1}{r} \Lb Z^k \phi \cdot L Z^k\phi .
\end{split}
\end{equation*}
We multiply both sides of the equation by $\delta^{2l}$ to normalize the contribution from the initial data to be close to $1$, therefore, we obtain
\begin{equation}\label{E 3 Lb}
\begin{split}
\delta^{2l}\int_{C_{u}}|\nablaslash Z^k \phi|^{2}&+\delta^{2l}\big|\int_{\underline{C}_{\underline{u}}}|\Lb Z^k \phi|^{2} \lesssim I_3^2 + \delta^{2l}\big|\doubleint_{D_{\ub,u}}Q(\nabla Z^k \phi, \nabla \phi)\Lb Z^k \phi \big| \\
&  + \sum_{p+q \leq k,\atop p<k, q<k} \delta^{2l}\big|\doubleint_{D_{\ub,u}}
Q(\nabla Z^{p}\phi, \nabla Z^{q}\phi)\Lb Z^i \phi \big|+\delta^{2l}\big|\doubleint_{D_{\ub,u}} \frac{1}{r} \Lb Z^k \phi \cdot L Z^k\phi\big|.
\end{split}
\end{equation}
We rewrite the right-hand side of the above inequality as
\begin{equation*}
I_3^2 +  S + T +W.
\end{equation*}
where $S$, $T$ and $W$ denote the three bulk integral terms in \eqref{E 3 Lb}. We will bound $S$, $T$ and $W$ one by one.

\bigskip

We begin with $S$, by definition, $S$ is bounded by the sum of the following integrals:
\begin{equation*}
\begin{split}
S_1 &=\delta^{2l} \doubleint_{{D_{\ub,u}}}\big(|L\phi|+|\nablaslash \phi|\big)|\Lb Z^k \phi|^2, \\
S_2 &=\delta^{2l} \doubleint_{{D_{\ub,u}}}|\Lb\phi||LZ^k\phi||\Lb Z^k \phi|,\\
S_3 &=\delta^{2l} \doubleint_{{D_{\ub,u}}}|\Lb \phi||\nablaslash Z^k\phi||\Lb Z^k \phi|, \\
S_4 &=\delta^{2l} \doubleint_{{D_{\ub,u}}}|\nablaslash\phi||L Z^k \phi||\Lb Z^k \phi|,\\
S_5 &=\delta^{2l} \doubleint_{{D_{\ub,u}}}(|L\phi|+|\nablaslash \phi|)|\nablaslash Z^k\phi||\Lb Z^k \phi|.
\end{split}
\end{equation*}
It suffices to bound the $S_i$'s one by one.

For $S_1$, in view of the $L^\infty(S_{\ub,u})$ estimates on $L\phi$ and $\nablaslash \phi$, we have
\begin{align*}
S_1 &\leq \delta^{2l}\int_{1-u}^{\ub}\int_{0}^{u}(\|L\phi\|_{L^{\infty}(S_{\ub',u'})}+ \|\nablaslash\phi\|_{L^{\infty}(S_{\ub',u'})}\|\Lb{Z^k}\phi\|_{L^{2}(S_{\ub',u'})}^{2} d u' d\ub'\\
&\leq \delta^{2l}\int_{1-u}^{\ub}\int_{0}^{u}\delta^{\frac{1}{4}}|\ub'|^{-\frac{6+\alpha}{4}}M\|\Lb{Z^k}\phi\|_{L^{2}(S_{\ub',u'})}^{2} d u'd\ub'\\
 &\lesssim \int_{1-u}^{u} \delta^{\frac{1}{4}}|\ub'|^{-\frac{6+\alpha}{4}}M  \big(\delta^{2l}\|\Lb{Z^k}\phi\|_{L^{2}(\Cb_{\ub'})}^{2}\big) d\ub'.
\end{align*}
In view of the bootstrap assumption, we bound $\delta^{2l}\|\Lb{Z^k}\phi\|_{L^{2}(\Cb_{\ub'})}^{2}$ by $M^2$. After an integration over $\ub'$ on $[1-u,\ub]$, we have
\begin{align*}
S_1 \lesssim \delta^{\frac{1}{4}}(|1-u|^{-\frac{2+\alpha}{4}}-|\ub|^{-\frac{2+\alpha}{4}})M^3.
\end{align*}
Because of $u\in [0,\delta]$, for sufficiently small $\delta$, we have
\begin{equation}\label{S1 Lb}
S_1 \lesssim \delta^{\frac{1}{4}}M^3.
\end{equation}

For $S_2$, we have
\begin{equation*}
S_{2}\lesssim \delta^{2l}\int_{1-u}^\ub \|L {Z^k}\phi\|_{L^{\infty}_{u}L^{2}(S_{\ub',u})} \|\Lb \phi \|_{L^{2}_{u}L^{\infty}(S_{\ub',u})}\|\Lb {Z^k}\phi\|_{L^{2}(\Cb_{\ub'})}d\ub'.
\end{equation*}
According to the bootstrap assumption \eqref{bootstrap assumption} and the estimates \eqref{preliminary estimates}, we have $\|L {Z^k}\phi\|_{L^{\infty}_{u}L^{2}(S_{\ub',u})} \lesssim 	\delta^{\frac{1}{2}-l}|\ub|^{-1}M$, $\|\Lb \phi \|_{L^{2}_{u}L^{\infty}(S_{\ub',u})}\lesssim  I_3|\ub|^{-1}$ and $\|\Lb {Z^k}\phi\|_{L^{2}(\Cb_{\ub'})}\lesssim \delta^{-l}M$, therefore, we can conclude that
\begin{equation}\label{S2 Lb}
S_2 \lesssim \delta^{\frac{1}{2}}M^2.
\end{equation}

For $S_3$, we have
\begin{align*}
S_{3}&\lesssim \delta^{2l}\int_{1-u}^\ub \|\nablaslash {Z^k}\phi\|_{L^{\infty}_{u}L^{2}(S_{\ub',u})} \|\Lb \phi \|_{L^{2}_{u}L^{\infty}(S_{\ub',u})}\|\Lb {Z^k}\phi\|_{L^{2}(\Cb_{\ub'})}d\ub'\\
&\lesssim I_{3}\cdot M \cdot \delta^{l} \int_{1-u}^\ub \|\nablaslash {Z^k}\phi\|_{L^{\infty}_{u}L^{2}(S_{\ub',u})} |\ub|^{-1} d\ub'.
\end{align*}
According to the $L^4$ estimates on $\nablaslash {Z^k}\phi$ on $S_{\ub,u}$ in \eqref{preliminary estimates}, we have $\|\nablaslash {Z^k}\phi\|_{L^2(S_{\ub,u})} \lesssim \delta^{\frac{1}{4}-l}|\ub|^{-\frac{1}{2}-\alpha}M$, this leads to
\begin{equation}\label{S3 Lb}
S_3 \lesssim \delta^{\frac{1}{4}}M^2.
\end{equation}

For $S_4$, we can proceed exactly as for $S_2$ (we just replace the factor $\Lb\phi$ in $S_2$ by $\nablaslash \phi$), this gives
\begin{equation}\label{S4 Lb}
S_4 \lesssim \delta^{\frac{5}{4}}M^3.
\end{equation}

For $S_5$, we can proceed exactly as for $S_3$ (we just replace the factor $|\Lb\phi|$ in $S_2$ by $|\nablaslash \phi|+|L\phi|$), this gives
\begin{equation}\label{S5 Lb}
S_5 \lesssim \delta M^3.
\end{equation}

\bigskip

We now estimate the second term in \eqref{E 3 Lb}, i.e. the estimates on $T$. According to the structure of null forms, we have
\begin{align*}
 T&=\sum_{p+q \leq k,\atop p<k, q<k} \delta^{2l}\big|\doubleint_{D_{\ub,u}} Q(\nabla Z^{p}\phi, \nabla Z^{q}\phi)\Lb Z^k \phi \big|\\
 &\lesssim \sum_{p+q \leq k,\atop p<k, q<k} \delta^{2l}\doubleint_{D_{\ub,u}} |\partial Z^{p}\phi| |\partial_g  Z^{q}\phi||\Lb Z^k \phi|,
\end{align*}
where $\partial\in\{\nablaslash,\Lb\}$ and  $\partial_g \in\{\nablaslash, L\}$. For each given term in the above summ, let $l'$ and $l''$ be total numbers of bad commutator $Z_b$'s appearing in $Z^p$ and $Z^q$ respectively. We remark that $l' + l'' \leq l$. Since $q\leq 1$, we have
\begin{align*}
T & \lesssim \delta^{2l}\sum_{p+q \leq k,\atop p<k, q<k}\int_{1-u}^{\ub}\int_{0}^{u}\|\partial Z^{p}\phi\|_{L^{4}(S_{\ub',u'})}\|\partial_{g}Z^{q}\phi\|_{L^{4}(S_{\ub',u'})}\|\Lb Z^{k}\phi\|_{L^{2}(S_{\ub',u'})}du'd\ub'
\end{align*}
By the second of \eqref{L2 Sobolev}, we have:
\begin{align}\label{partialg L4}
\|\partial_{g}Z^{q}\phi\|_{L^{4}(S_{u',\ub'})}\lesssim\ub^{\prime-1/2}\|\Lb \partial_{g}Z^{q}\phi\|_{L^{2}(\Cb_{\ub})}^{1/2}\left(\|\partial_{g}Z^{q}\phi\|_{L^{2}(\Cb_{\ub})}+\|\Omega\partial_{g}Z^{q}\phi\|_{L^{2}(\Cb_{\ub})}\right)^{1/2}
\end{align}
Now by \eqref{formulas for partial}, we have schematically:
\begin{align*}
L\sim\frac{1}{\ub}S+\frac{1}{\ub}\frac{x^{i}}{r}\Omega_{0i},\quad \slashed{\nabla}\sim\frac{1}{\ub}\Omega_{ij}
\end{align*}
which imply, for any smooth function $f$:
\begin{align*}
\Lb\partial_{g}f\sim \frac{1}{\ub}\Lb Z_{g}f,
\quad \Omega f\sim\ub\slashed{\nabla}f,\quad \partial_{g}f\sim\frac{1}{\ub}Z_{g}f
\end{align*}
If $\partial_{g}=\slashed{\nabla}$, then the second factor on the right hand side of \eqref{partialg L4} is bounded through the bootstrap assumption \eqref{bootstrap assumption} by:
\begin{align*}
\left(\|\slashed{\nabla}Z^{q}\phi\|_{L^{2}(\Cb_{\ub})}+\|\slashed{\nabla}Z_{g}Z^{q}\phi\|_{L^{2}(\Cb_{\ub})}\right)^{1/2}\lesssim \delta^{1/4-l''/2}M^{1/2}
\end{align*}
If $\partial_{g}=L$, by virtue of \eqref{L-2 for L Zk phi} the first term in the parenthesis is bounded by:
\begin{align*}
\|LZ^{q}\phi\|_{L^{2}(\Cb_{\ub})}\lesssim\delta^{1-l''}\ub^{-1}M
\end{align*}
while for the second term we have:
\begin{align*}
\|\Omega\partial_{g}Z^{q}\phi\|_{L^{2}(\Cb_{\ub})}\lesssim\|\slashed{\nabla}Z_{g}Z^{q}\phi\|_{L^{2}(\Cb_{\ub})}\lesssim\delta^{1/2-l''}M^{1/2}
\end{align*}
These together with \eqref{partialg L4} imply:
\begin{align}\label{partialg-L4 estimate}
\|\partial_{g}Z^{q}\phi\|_{L^{4}(S_{u',\ub'})}\lesssim \delta^{1/4-l''}M\ub^{\prime-1}
\end{align}
Therefore \eqref{L4 Sobolev} implies:
\begin{align*}
T&\lesssim \delta^{1/4}M\sum_{p+q \leq k,\atop p<k, q<k}\int_{1-u}^{\ub}\ub^{\prime-3/2}\int_{0}^{u}\left(\delta^{l-l''}\|\partial Z^{q}Z_{g}^{q'}\phi\|_{L^{2}(S_{u'\ub'})}\delta^{l}\|\Lb Z^{k}\phi\|_{L^{2}(S_{u'\ub'})}\right)du'd\ub'\\
&\lesssim\delta^{1/4}M\sum_{p+q \leq k,\atop p<k, q<k}\int_{1-u}^{\ub}\ub^{\prime-3/2}\left(\delta^{l-l''}\|\partial Z^{q}Z_{g}^{q'}\phi\|_{L^{2}(\Cb_{\ub'})}\delta^{l}\|\Lb Z^{k}\phi\|_{L^{2}(\Cb_{\ub'})}\right)d\ub'\\
&\lesssim \delta^{1/4}M^{3}\int_{1-u}^{\ub}\ub^{\prime-3/2}d\ub'
\end{align*}
where $|q'|\leq 1$.
This eventually yields
\begin{equation}\label{T Lb}
T \lesssim \delta^{\frac{1}{4}}M^3.
\end{equation}

\bigskip

It remains to bound the third term $W$ in \eqref{E 3 Lb}. It is similar to $T_2$. We simply bound  $L Z^k \phi $ and $\Lb Z^k \phi$ on $\Cb_{\ub}$. Although the bound of $|L Z^{p}\phi|$ on $\Cb_{\ub}$ is not directly from the bootstrap assumption, in view of \eqref{preliminary estimates} and the fact that $k\leq 2$, we can bound $L Z^{k}\phi$ first on $L^2(S_{\ub,u})$ and then on $L^2(\Cb_{\ub})$. This leads to
\begin{equation}\label{W Lb}
 W \lesssim  \delta M^2.
\end{equation}

By combining \eqref{E 3 Lb} with \eqref{S1 Lb}, \eqref{S2 Lb}, \eqref{S3 Lb}, \eqref{S4 Lb}, \eqref{S5 Lb}, \eqref{T Lb} and \eqref{W Lb}, for sufficiently small $\delta$, we obtain
\begin{equation*}
\delta^{2l}\int_{C_{u}}|\nablaslash Z^k \phi|^{2} +\delta^{2l}\int_{\underline{C}_{\underline{u}}}|\Lb Z^k \phi|^{2} \lesssim I_3^2 + \delta^{\frac{1}{4}}M^3.
\end{equation*}
In other words, for all $0\leq l \leq k \leq 2$, we have
\begin{equation}\label{Eenergy estimates 3 Lb}
\delta^{l}\|\nablaslash Z_b^l Z_g^{k-l} \phi\|_{L^2(C_u)} +\delta^{l}\|\Lb Z_b^l Z_g^{k-l} \phi\|_{L^2(\Cb_{\ub})} \lesssim I_3 + \delta^{\frac{1}{8}}M^{\frac{3}{2}}.
\end{equation}

In the fundamental energy identity \eqref{fundamental energy identity} and \eqref{commuated main equatin for energy estimates}, we replace $\phi$ by $Z^k \phi$ and take $X = \ub^\alpha L$ to obtain
\begin{equation*}
\begin{split}
\int_{C_{u}}|\ub|^{\alpha}&|L {Z^k}{\phi}|^{2}+\int_{\underline{C}_{\underline{u}}}|\ub|^{\alpha}|\nablaslash {Z^k} {\phi}|^{2}= \int_{\Sigma_1}|	\ub|^\alpha \big(|\nablaslash Z^k \phi|^{2}+|L Z^k \phi|^{2}\big) +\doubleint_{\mathcal{D}_{\ub,u}}|\ub|^{\alpha} Q(\nabla Z^k \phi, \nabla\phi)L {Z^k} {\phi} \\
& + \sum_{p+q \leq k,\atop p<k, q<k} \doubleint_{D_{\ub,u}} |\ub|^{\alpha} Q(\nabla Z^{p}\phi, \nabla Z^{q}\phi)L {Z^k}{\phi} +\doubleint_{D_{\ub,u}} \frac{|\ub|^{\alpha}}{r} \Lb {Z^k} {\phi} \cdot L {Z^k}{\phi}\\
& -2\alpha \doubleint_{D_{\ub,u}}|\ub|^{\alpha-1}|\nablaslash {Z^k} {\phi}|^{2}.
\end{split}
\end{equation*}
We multiply both sides of the equation by $\delta^{2l-1}$ to renormalize the contribution from the initial data to be close to $1$. We remark that this normalization is respect to the relaxed estimates on $\nablaslash \phi$. By dropping of the last 
negative term in the above equation, we obtain
\begin{equation}\label{E 3 L}
\begin{split}
\delta^{2l-1}\int_{C_{u}}|\ub|^{\alpha}&|L {Z^k}{\phi}|^{2}+\delta^{2l-1} \int_{\underline{C}_{\underline{u}}}|\ub|^{\alpha}|\nablaslash {Z^k} {\phi}|^{2}\lesssim I_3^2  +\delta^{2l-1}\big|\doubleint_{\mathcal{D}_{\ub,u}}|\ub|^{\alpha} Q(\nabla Z^k \phi, \nabla\phi)L {Z^k} {\phi}\big| \\
& + \sum_{p+q \leq k,\atop p<k, q<k} \delta^{2l-1} \big|\doubleint_{D_{\ub,u}} |\ub|^{\alpha} Q(\nabla Z^{p}\phi, \nabla Z^{q}\phi)L {Z^k}{\phi}\big| +\delta^{2l-1}\big|\doubleint_{D_{\ub,u}} \frac{|\ub|^{\alpha}}{r} \Lb {Z^k} {\phi} \cdot L {Z^k}{\phi}\big|,
\end{split}
\end{equation}
We rewrite the right-hand side of the above inequality as
\begin{equation*}
I_3^2 +  S + T +W.
\end{equation*}
where $S$, $T$ and $W$ denote the three bulk integral terms in \eqref{E 3 L}. We now bound $S$, $T$ and $W$ one by one.

\bigskip

We begin with $S$. According to the definition of $S$ and the structure \eqref{null form bound} for null forms, $S$ is bounded by the sum of the the following terms:
\begin{equation*}
\begin{split}
 S_1 &= \delta^{2l-1}\doubleint_{\mathcal{D}_{\ub,u}}|\ub|^{\alpha}\big(|\Lb {\phi}|+|\nablaslash\phi|\big)|L{Z^k} {\phi}|^2, \\
 S_2 &=\delta^{2l-1} \doubleint_{\mathcal{D}_{\ub,u}}|\ub|^{\alpha}\big(|\nablaslash{\phi}|+|L{\phi}|\big)|\Lb{Z^k}{\phi}||L{Z^k} {\phi}|,\\
 S_3 &= \delta^{2l-1}\doubleint_{\mathcal{D}_{\ub,u}}|\ub|^{\alpha}\big(|L {\phi}|+|\nablaslash {\phi}|\big)|\nablaslash{Z^k}{\phi}||L{Z^k} {\phi}|, \\
 S_4 &= \delta^{2l-1}\doubleint_{\mathcal{D}_{\ub,u}}|\ub|^{\alpha}|\Lb{\phi}||\nablaslash{Z^k}{\phi}||L{Z^k} {\phi}|.
\end{split}
\end{equation*}

The idea to bound the $S_i$'s are exactly the same as before. Roughly speaking, we bound all the first order derivative components of $\nabla \phi$ in $L^\infty(S_{\ub,u})$.

For $S_1$,  we have
\begin{align*}
 S_1 &\leq \delta^{2l-1} \doubleint_{\mathcal{D}_{\ub,u}}|\ub|^{\alpha} \big( |\ub|^{-1} \delta^{-\frac{1}{2}}M\big) |L{Z^k} {\phi}|^2 \\
 &\lesssim \delta^{-\frac{1}{2}}M\int_{0}^u \Big(\delta^{2l-1}\int_{C_{u'}} |\ub|^\alpha |L{Z^k} {\phi}|^2\Big)du'.
\end{align*}
According to the bootstrap assumption on $\delta^{l-\frac{1}{2}} \|L Z^k \phi\|_{L^2(C_u)}$, we obtain
\begin{equation}\label{S1 L}
S_1 \lesssim \delta^{\frac{1}{2}}M^3.
\end{equation}

For $S_{2}$, since $k \leq 2$, we use the bound on $L Z^k \phi$ on $S_{\ub,u}$ to derive 
$\delta^{l}\|L Z^k \phi\|_{L^2(\Cb_{\ub'})} \lesssim \delta |\ub|^{-1} M$. Therefore, we can proceed as follows:
\begin{align*}
 S_2 &\leq \delta^{2l-1}\int_{1-u}^{\ub}|\ub|^{\alpha} \big( |\ub|^{-\frac{3}{2}-\frac{\alpha}{4}} \delta^{\frac{1}{4}}M\big) \|\Lb Z^k \phi\|_{L^2(\Cb_{\ub'})}\|L Z^k \phi\|_{L^2(\Cb_{\ub'})} d\ub.
\end{align*}
According to the bootstrap assumptions, we finally obtain
\begin{equation}\label{S2 L}
S_2 \lesssim \delta^{\frac{1}{4}}M^3.
\end{equation}

The estimates on $S_3$ can be obtained in a similar way as $S_2$: we simply replace $\Lb Z^k \phi$ by $\nablaslash Z^k \phi$ and proceed exactly the same as before. This gives
\begin{equation}\label{S3 L}
S_3 \lesssim \delta^{\frac{3}{4}}M^3.
\end{equation}

For $S_4$, we first make the following remark:

\begin{remark}
It seems to be natural to derive the estimates by putting $\nablaslash Z^k \phi$ in the ${L^2(\Cb_{\ub'})}$ norm. In fact, this does not work due to the fact that we have relaxed the estimates on the rotational directions. To illustrate the idea, we may proceed as follows:
\begin{align*}
 S_4 &\leq \delta^{2l-1}\int_{1-u}^{\ub}|\ub|^{\alpha} \big( |\ub|^{-1} \delta^{-\frac{1}{2}}I_3\big) \|\nablaslash Z^k \phi\|_{L^2(\Cb_{\ub'})}\|L Z^k \phi\|_{L^2(\Cb_{\ub'})} d\ub\\
 &\leq \delta^{-1}\int_{1-u}^{\ub}|\ub|^{\frac{\alpha}{2}} \big( |\ub|^{-1} \delta^{-\frac{1}{2}}I_3\big)\big( \delta^{\frac{1}{2}}M \big)\big(\delta |\ub|^{-1} M\big) d\ub'\\
 &\lesssim M^2.
\end{align*}
This estimate is certainly not good since we do not have a $\delta$ (to some positive power) factor in front of the possibly large constant $M$.

At this point, we have to use the bootstrap assumptions on the fourth order derivatives of $\phi$ to improve the relaxed estimates on $\nablaslash$-direction.
\end{remark}

The above remark suggests to put $\nablaslash Z^k {\phi}$ in $L^4(S_{\ub,u})$ norm to get an extra $\delta^\frac{1}{4}$ factor. In fact, we have
\begin{align*}
 S_4 &\leq \delta^{2l-1}\int_{1-u}^{\ub}|\ub|^{\alpha} \|\Lb{\phi}\|_{L^{2}_{u}L^{4}(S_{\ub',u})} \|\nablaslash Z^k {\phi}\|_{L^{\infty}_{u}L^{4}(S_{\ub',u}))} \|L Z^k \phi\|_{L^2(\Cb_{\ub'})}d \ub'.
\end{align*}
Since we have already derived estimates on 
$\|\Lb \phi\|_{L^{4}(S_{\ub,u})}$ and $\|L Z^k \phi\|_{L^2(S_{\ub,u})}$ $(k\leq 2)$, a direct computation yields
\begin{equation}\label{S4 L}
S_4 \lesssim \delta^{\frac{1}{4}}M^2.
\end{equation}

\bigskip

We turn to the estimates on $T$. According to the structure of null forms, we have
\begin{align*}
 T &\lesssim \delta^{2l-1}\sum_{p+q \leq k,\atop p<k, q<k} \doubleint_{D_{\ub,u}} |\ub|^\alpha |\partial Z^{p}\phi| |\partial_g  Z^{q}\phi||L Z^k \phi|,
\end{align*}
where $\partial\in\{\nablaslash, \Lb\}$ and  $\partial_g \in\{\nablaslash, L\}$.

Here we postpone the estimate for $T$ until we estimate the top order energy $E_{\leq3}(\ub,u)$, because the estimates for $T$ corresponding to $E_{\leq2}$ and $E_{\leq3}$ are identical. Instead, we just state the result:
\begin{equation}\label{T L}
T \lesssim \delta^{-1}\sum_{p+q\leq k,\atop p<k,q<k}\int_{0}^{u}\delta^{2l-1}\|\ub^{\prime\alpha/2}LZ^{k}\phi\|^{2}_{L^{2}(C_{u'})}du'+I_{3}^{4}
\end{equation}

\bigskip

It remains to control $W = \delta^{2l-1}\int\!\!\!\int_{\mathcal{D}_{\ub,u}} \frac{|\ub|^{\alpha}}{r}|L{Z^k} {\phi}||\Lb{Z^k}{\phi}|$. We proceed as follows
\begin{align*}
W &\lesssim \delta^{2l-1} \doubleint_{\mathcal{D}_{\ub,u}}\big(\delta^{-\frac{1}{2}}|\ub|^{\frac{\alpha}{2}} |L{Z^k}{\phi}|\big) \big(\delta^{\frac{1}{2}} |\ub|^{-\frac{2-\alpha}{2}} |\Lb{Z^k}{\phi}|\big) \\
&\stackrel{Cauchy-Schwarz}{\lesssim} \delta^{2l-1}\Big( \doubleint_{\mathcal{D}_{\ub,u}}\frac{\delta}{|\ub'|^{2-\alpha}}|\Lb{Z^k}{\phi}|^2+  \doubleint_{\mathcal{D}_{\ub,u}} \frac{1}{\delta}||\ub'|^{\frac{\alpha}{2}}L{Z^k}{\phi}|^2\Big)\\
&= \delta^{2l}\int_{1-u}^{\ub}\frac{1}{|\ub'|^{2-\alpha}} \|\Lb{Z^k}{\phi}\|^2_{L^{2}(C_{u'})} d\ub' + \delta^{2l-2}\int_{0}^{u}\||\ub|^{\frac{\alpha}{2}}L{Z^k}{\phi}\|^2_{L^2(C_{u'})} d u'.
\end{align*}
The first term in the last line has already been controlled in \eqref{Eenergy estimates 3 Lb}. In view of the fact that $\alpha <1$ (this is crucial to make the first factor integrable in $\ub$!), for sufficiently small $\delta$, we obtain
\begin{equation}\label{W L}
W \lesssim I_{3}^2 +  \delta^{-1} \int_{0}^{u} \delta^{2l-1} \||\ub|^{\frac{\alpha}{2}}L{Z^k}{\phi}\|^2_{L^2(C_{u'})} d u'.
\end{equation}

\bigskip

By combining \eqref{E 3 L} with \eqref{S1 L}, \eqref{S2 L}, \eqref{S3 L}, \eqref{S4 L}, \eqref{T L} and \eqref{W L}, for sufficiently small $\delta$, we obtain
\begin{equation*}
\delta^{2l-1}\int_{C_{u}}|\ub|^{\alpha} |L {Z^k}{\phi}|^{2}+\delta^{2l-1} \int_{\underline{C}_{\underline{u}}}|\ub|^{\alpha}|\nablaslash {Z^k} {\phi}|^{2}\lesssim I_3^4 + \delta^{\frac{1}{4}}M^3 + \delta^{-1} \int_{0}^{u} \delta^{2l-1} \||\ub|^{\frac{\alpha}{2}}L{Z^k}{\phi}\|^2_{L^2(C_{u'})} d u'
\end{equation*}
The last term on the right-hand side can be removed by the Gronwall's inequality. This finally proves that, for all $0\leq l \leq k \leq 2$, we have
\begin{equation}\label{Eenergy estimates 3 L}
\delta^{l-\frac{1}{2}}\||\ub|^\frac{\alpha}{2}L Z_b^l Z_g^{k-l} \phi\|_{L^2(C_u)} +\delta^{l-\frac{1}{2}}\|\nablaslash Z_b^l Z_g^{k-l} \phi\|_{L^2(\Cb_{\ub})} \lesssim C(I_3) + \delta^{\frac{1}{8}}M^{\frac{3}{2}}.
\end{equation}

\bigskip

The estimates \eqref{Eenergy estimates 3 Lb} and \eqref{Eenergy estimates 3 L} together implies
\begin{equation}\label{energy estimates k leq 2}
E_{\leq 2}(u,\ub) + \Eb_{\leq 2}( u,\ub) \leq C(I_3)+\delta^{\frac{1}{8}}M^{\frac{3}{2}}.	
\end{equation}

\subsubsection{Estimates on $E_{3}$ and $\Eb_{3}$}\label{third order energy estimates}
We take $k=3$ in \eqref{commuated main equatin for energy estimates}.  Let $l\leq k$ be the number of $Z_b$'s appearing in $Z^3$, i.e. $Z^3 = Z_b^l Z_g^{3-l}$. We take $Z^3 \phi$ in the place of $\phi$ in \eqref{fundamental energy identity} and  take the multiplier $X = \Lb$, this yields
\begin{equation*}
\begin{split}
\int_{C_{u}}|\nablaslash Z^3 \phi|^{2}&+\int_{\underline{C}_{\underline{u}}}|\Lb Z^3 \phi|^{2} =\int_{\Sigma_1}|\nablaslash Z^3 \phi|^{2}+|\Lb Z^3 \phi|^{2}  + \doubleint_{D_{\ub,u}}Q(\nabla Z^3 \phi, \nabla \phi)\Lb Z^3 \phi  \\
&  + \sum_{p+q \leq 3,\atop p<3, q<3} \doubleint_{D_{\ub,u}}
Q(\nabla Z^{p}\phi, \nabla Z^{q}\phi)\Lb Z^3 \phi -\doubleint_{D_{\ub,u}} \frac{1}{r} \Lb Z^3 \phi \cdot L Z^3\phi .
\end{split}
\end{equation*}
After a renormalization in $\delta$, we obtain
\begin{equation}\label{E 3 Lbb}
\begin{split}
\delta^{2l}\int_{C_{u}}|\nablaslash Z^3 \phi|^{2}&+\delta^{2l}\int_{\underline{C}_{\underline{u}}}|\Lb Z^3 \phi|^{2} \lesssim I_4^2 + \delta^{2l}\big|\doubleint_{D_{\ub,u}}Q(\nabla Z^3 \phi, \nabla \phi)\Lb Z^3 \phi \big| \\
&  + \sum_{p+q \leq 3,\atop p<3, q<3} \delta^{2l}\big|\doubleint_{D_{\ub,u}}
Q(\nabla Z^{p}\phi, \nabla Z^{q}\phi)\Lb Z^3 \phi \big|+\delta^{2l}\big|\doubleint_{D_{\ub,u}} \frac{1}{r} \Lb Z^3 \phi \cdot L Z^3\phi\big|\\
&= I_4^2 +  S + T +W,
\end{split}
\end{equation}
where $S$, $T$ and $W$ denote the three bulk integral terms. We will bound $S$, $T$ and $W$ one by one.

\bigskip

We start with $S$. It  can be bounded by the sum of the following terms:
\begin{align*}
S_1 &= \delta^{2l}\doubleint_{\mathcal{D}_{\ub,u}}|\Lb {\phi}|\big(|L Z^3 {\phi}|+|\nablaslash Z^3{\phi}|\big)|\Lb Z^3 {\phi}|, \\
S_2 &= \delta^{2l}\doubleint_{\mathcal{D}_{\ub,u}}|L{\phi}|\big(|\nablaslash Z^3 {\phi}|+|\Lb Z^3 {\phi}|\big)|\Lb Z^3 {\phi}|,\\
S_3 &= \delta^{2l}\doubleint_{\mathcal{D}_{\ub,u}}|\nablaslash {\phi}|\big(|L Z^3 {\phi}|+|\Lb Z^3 {\phi}|+|\nablaslash Z^3 {\phi}|\big)|\Lb Z^3{\phi}|.
\end{align*}

For $S_{1}$,  according to the $L^\infty(S_{\ub,u})$ estimates on  $\Lb{\phi}$, we have
\begin{align*}
S_{1}&\lesssim \delta^{2l-\frac{1}{2}}\doubleint_{\mathcal{D}_{\ub,u}} |\ub|^{-1}|\nablaslash Z^{3}{\phi}| |\Lb Z^{3}{\phi}|+ \delta^{2l-\frac{1}{2}}\doubleint_{\mathcal{D}_{\ub,u}} |\ub|^{-\frac{1+\alpha}{2}} \cdot \big(|\ub|^{\frac{\alpha}{2}} L Z^{3}{\phi}|\big)\cdot\big(|\ub|^{-\frac{1}{2}}|\Lb Z^{3}{\phi}|\big)\\
&=S_{11}+S_{12}.
\end{align*}
For $S_{11}$, according to Cauchy-Schwarz inequality, we have
\begin{equation}\label{S11}
\begin{split}
S_{11}& \lesssim \delta^{2l-1}\doubleint_{\mathcal{D}_{\ub,u}}|\nablaslash Z^{3}{\phi}|^{2}+\delta^{2l}\doubleint_{\mathcal{D}_{\ub,u}}|\ub|^{-2}|\Lb Z^{3}{\phi}|^{2}\\
&\lesssim \int_{0}^{u}\frac{1}{\delta} \cdot \delta^{2l}\int_{C_{u'}}|\nablaslash Z^{3}{\phi}|^{2} du'+\int_{1-u}^{\ub} |\ub'|^{-2} \delta^{2l}\int_{\Cb_{\ub}}|\Lb Z^{3}{\phi}|^{2} d\ub'.
\end{split}
\end{equation}
For $S_{12}$, we still use Cauchy-Schwarz inequality to derive
\begin{align*}
S_{12}&\lesssim \delta^{2l-1}\doubleint_{\mathcal{D}_{\ub,u}}\big(|\ub|^{\frac{\alpha}{2}}|L Z^{3}{\phi}|\big)^{2}+\delta^{2l} \doubleint_{\mathcal{D}_{\ub,u}}|\ub|^{-2-\alpha} |\Lb Z^{3}{\phi}|^{2}\\
&\lesssim \int_{0}^{u}\frac{1}{\delta} \cdot \delta^{2l}\int_{C_{u'}}|LZ^{3}{\phi}|^{2} du'+\int_{1-u}^{\ub} |\ub'|^{-2-\alpha} \delta^{2l}\int_{\Cb_{\ub}}|\Lb Z^{3}{\phi}|^{2} d\ub'\\
&\lesssim \delta M^{2}+\int_{1-u}^{\ub} |\ub'|^{-2-\alpha} \delta^{2l}\int_{\Cb_{\ub}}|\Lb Z^{3}{\phi}|^{2} d\ub'
\end{align*}
Therefore, we obtain
\begin{equation}\label{3rd order S1 L}
S_1 \lesssim \delta M^2 + \int_{0}^{u}\frac{1}{\delta} \cdot \delta^{2l}\int_{C_{u'}}|\nablaslash Z^{3}{\phi}|^{2} du'+\int_{1-u}^{\ub} |\ub'|^{-2} \delta^{2l}\int_{\Cb_{\ub}}|\Lb Z^{3}{\phi}|^{2} d\ub'.
\end{equation}

For $S_{2}$,  according to the $L^\infty(S_{\ub,u})$ estimates on  $L{\phi}$, we have
\begin{align*}
S_{2}&\lesssim \delta^{2l+\frac{1}{2}}M\doubleint_{\mathcal{D}_{\ub,u}} |\ub|^{-2}|\nablaslash Z^{3}{\phi}| |\Lb Z^{3}{\phi}|+ \delta^{2l+\frac{1}{2}}M\doubleint_{\mathcal{D}_{\ub,u}} |\ub|^{-2}|\Lb Z^{3}{\phi}|^2\\
&=S_{21}+S_{22}.
\end{align*}
For $S_{21}$, since $1 \lesssim |\ub|$, we have
\begin{align*}
S_{21}& \lesssim \delta^{2l}M\doubleint_{\mathcal{D}_{\ub,u}}|\nablaslash Z^{3}{\phi}|^{2}+\delta^{2l+1}M\doubleint_{\mathcal{D}_{\ub,u}}|\ub|^{-2}|\Lb Z^{3}{\phi}|^{2}\\
&\lesssim \int_{0}^{u}\delta^{2l}M\int_{C_{u'}}|\nablaslash Z^{3}{\phi}|^{2} du'+\delta M\int_{1-u}^{\ub} |\ub'|^{-2} \delta^{2l}\int_{\Cb_{\ub}}|\Lb Z^{3}{\phi}|^{2} d\ub'\\
&\lesssim \delta M^3.
\end{align*}
For $S_{22}$, we have
\begin{align*}
S_{22}&\lesssim \delta^{\frac{1}{2}} M \int_{1-u}^{\ub} |\ub'|^{-2} \delta^{2l}\int_{\Cb_{\ub'}}|\Lb Z^{3}{\phi}|^{2} d\ub'\\
&\lesssim \delta^{\frac{1}{2}}M^3.
\end{align*}
Therefore, we obtain
\begin{equation}\label{3rd order S2 L}
S_2 \lesssim  \delta^{\frac{1}{2}}M^3.
\end{equation}

For $S_{3}$,  according to the $L^\infty(S_{\ub,u})$ estimates on  $\nablaslash {\phi}$, it is bounded by the following three terms:
\begin{align*}
S_{31} &= \delta^{2l+\frac{1}{4}}M\doubleint_{\mathcal{D}_{\ub,u}} |\ub|^{-\frac{3}{2}-\frac{\alpha}{4}}|L Z^{3}{\phi}| |\Lb Z^{3}{\phi}|,\\
S_{32} &= \delta^{2l+\frac{1}{4}}M\doubleint_{\mathcal{D}_{\ub,u}} |\ub|^{-\frac{3}{2}-\frac{\alpha}{4}}|\Lb  Z^{3}{\phi}|^2,\\
S_{33} &=\delta^{2l+\frac{1}{4}}M\doubleint_{\mathcal{D}_{\ub,u}} |\ub|^{-\frac{3}{2}-\frac{\alpha}{4}}|\nablaslash Z^{3}{\phi}| |\Lb Z^{3}{\phi}|.
\end{align*}
To bound $S_{31}$, we follow exactly the same way for $S_{11}$, this yields
\begin{equation*}
S_{31} \lesssim  \delta^{\frac{1}{4}}M^3.
\end{equation*}
To bound $S_{32}$, we follow exactly the same way for $S_{22}$, this yields
\begin{equation*}
S_{32} \lesssim  \delta^{\frac{1}{4}}M^3.
\end{equation*}
To bound $S_{33}$, we follow exactly the same way for 
$S_{21}$, this yields
\begin{equation*}
S_{33} \lesssim  \delta^{\frac{1}{4}}M^3.
\end{equation*}
Therefore, we obtain
\begin{equation}\label{3rd order S3 L}
S_3 \lesssim  \delta^{\frac{1}{4}}M^3.
\end{equation}

\bigskip

We turn to the estimates on $T$. According to the structure of null forms, we have
\begin{align*}
 T &\lesssim \sum_{p+q \leq 3,\atop p\leq 2, q\leq 2} \delta^{2l}\doubleint_{D_{\ub,u}} |\partial Z^{p}\phi| |\partial_g  Z^{q}\phi||\Lb Z^3 \phi|.
\end{align*}
where $\partial\in\{\nablaslash, \Lb\}$ and  $\partial_g \in\{\nablaslash, L\}$.
By using exactly the same method as we derive \eqref{T Lb}, we obtain
\begin{equation}\label{3rd order T L}
T\lesssim 
\delta^\frac{1}{4}M^3 
\end{equation}

\bigskip

It remains to bound $W = \delta^{2l}\int\!\!\!\int_{\mathcal{D}_{\ub,u}} \frac{1}{r}|L{Z^3} {\phi}||\Lb{Z^3}{\phi}|$. According to Cauchy-Schwarz inequality, we have
\begin{align*}
W &{\lesssim} \delta^{2l}\Big( \doubleint_{\mathcal{D}_{\ub,u}}\frac{\delta}{|\ub'|^{2-\alpha}}|\Lb{Z^3}{\phi}|^2+  \doubleint_{\mathcal{D}_{\ub,u}} \frac{1}{\delta}||\ub'|^{\frac{\alpha}{2}}L{Z^3}{\phi}|^2\Big)\\
&=\int_{1-u}^{\ub}\frac{1}{|\ub'|^{2-\alpha}} \cdot \delta^{2l+1}\|\Lb{Z^3}{\phi}\|^2_{L^{2}(\Cb_{\ub'})} d\ub' + \delta^{2l-1}\int_{0}^{u}\||\ub|^{\frac{\alpha}{2}}L{Z^3}{\phi}\|^2_{L^2(C_{u'})} d u'.
\end{align*}
This yields $W \lesssim \delta M^2$. Combining this estimate with \eqref{3rd order S1 L}, \eqref{3rd order S2 L}, \eqref{3rd order S3 L} and \eqref{3rd order T L},
 we obtain
\begin{align*}
\delta^{2l}\int_{C_{u}}|\nablaslash Z^3 \phi|^{2}+\delta^{2l}\int_{\underline{C}_{\underline{u}}}|\Lb Z^3 \phi|&^{2} \lesssim I_4^2 + \delta^{\frac{1}{4}}M^3.
\end{align*}
In other words, for all $0\leq l \leq 3$, we have
\begin{equation}\label{Eenergy estimates top Lb}
\delta^{l}\|\nablaslash Z_b^l Z_g^{3-l} \phi\|_{L^2(C_u)} +\delta^{l}\|\Lb Z_b^l Z_g^{3-l} \phi\|_{L^2(\Cb_{\ub})} \lesssim I_4 + \delta^{\frac{1}{8}}M^{\frac{3}{2}}.
\end{equation}

\bigskip

Similar to the derivation for \eqref{E 3 L} (by taking $k=3$), we have
\begin{equation*}
\begin{split}
\delta^{2l-1}\int_{C_{u}}|\ub|^{\alpha}&|L {Z^3}{\phi}|^{2}+\delta^{2l-1} \int_{\underline{C}_{\underline{u}}}|\ub|^{\alpha}|\nablaslash {Z^3} {\phi}|^{2}\lesssim I_4^2  +\delta^{2l-1}\big|\doubleint_{\mathcal{D}_{\ub,u}}|\ub|^{\alpha} Q(\nabla Z^3 \phi, \nabla\phi)L {Z^3} {\phi}\big| \\
& + \sum_{p+q \leq 3,\atop p\leq 2, q \leq 2} \delta^{2l-1} \big|\doubleint_{D_{\ub,u}} |\ub|^{\alpha} Q(\nabla Z^{p}\phi, \nabla Z^{q}\phi)L {Z^3}{\phi}\big| +\delta^{2l-2}\big|\doubleint_{D_{\ub,u}} \frac{|\ub|^{\alpha}}{r} \Lb {Z^3} {\phi} \cdot L {Z^3}{\phi}\big|,
\end{split}
\end{equation*}
We rewrite the above inequality as
\begin{equation}\label{E top L}
\delta^{2l-1}\int_{C_{u}}|\ub|^{\alpha}|L {Z^3}{\phi}|^{2}+\delta^{2l-1} \int_{\underline{C}_{\underline{u}}}|\ub|^{\alpha}|\nablaslash {Z^3} {\phi}|^{2}\lesssim I_4^2  +  S + T +W.
\end{equation}
where $S$, $T$ and $W$ denote the three bulk integral terms in an obvious way. We now bound $S$, $T$ and $W$ one by one.

We begin with $S$ which is bounded by the sum of the the following terms:
\begin{equation*}
\begin{split}
 S_1 &= \delta^{2l-1}\doubleint_{\mathcal{D}_{\ub,u}}|\ub|^{\alpha}\big(|\Lb {\phi}|+|\nablaslash\phi|\big)|L{Z^3} {\phi}|^2, \\
 S_2 &=\delta^{2l-1} \doubleint_{\mathcal{D}_{\ub,u}}|\ub|^{\alpha}\big(|\nablaslash{\phi}|+|L{\phi}|\big)|\Lb{Z^3}{\phi}||L{Z^3} {\phi}|,\\
 S_3 &= \delta^{2l-1}\doubleint_{\mathcal{D}_{\ub,u}}|\ub|^{\alpha}\big(|L {\phi}|+|\nablaslash {\phi}|\big)|\nablaslash{Z^3}{\phi}||L{Z^3} {\phi}|, \\
 S_4 &= \delta^{2l-1}\doubleint_{\mathcal{D}_{\ub,u}}|\ub|^{\alpha}|\Lb{\phi}||\nablaslash{Z^3}{\phi}||L{Z^3} {\phi}|.
\end{split}
\end{equation*}

For $S_1$,  we have
\begin{align*}
 S_1 &\leq \delta^{2l-1} \doubleint_{\mathcal{D}_{\ub,u}}|\ub|^{\alpha} \big( |\ub|^{-1} \delta^{-\frac{1}{2}}M\big) |L{Z^3} {\phi}|^2 \\
 &\lesssim \delta^{-\frac{1}{2}}M\int_{0}^u \Big(\delta^{2l-1}\int_{C_{u'}} |\ub|^\alpha |L{Z^k} {\phi}|^2\Big)du'\\
 &\lesssim \delta^{\frac{1}{2}}M^3.
\end{align*}

For $S_{2}$, we use $L^\infty$ bound on $\nablaslash \phi$ and $L \phi$. Since we have already derived estimates on $E_{\leq 2}(u,\ub)$ and $\Eb_{\leq 2}(u,\ub)$, for sufficiently small $\delta$, we indeed have
\begin{equation}\label{improved L infty estimates}
|L {\phi}|+|\nablaslash\phi| \lesssim  |\ub|^{-\frac{3}{2}-\frac{\alpha}{4}}\delta^{\frac{1}{4}}C(I_3).
\end{equation}
We remark that this estimate is better than those in \eqref{preliminary estimates} since we have improved the big bootstrap constant $M$ to be a constant depending only on the size $I_{n}$ of the rescaled data. Therefore, according to Cauchy-Schwarz inequality, we have
\begin{align*}
 S_2 &\lesssim  C(I_3) \delta^{2l-1}\doubleint_{\mathcal{D}_{\ub,u}} \big( |\ub|^{-\frac{3}{2}-\frac{\alpha}{4}} \delta^{\frac{1}{4}}\big)|\Lb{Z^3}{\phi}||L{Z^3} {\phi}|\\
&\lesssim \delta^{\frac{1}{4}}\int_{1-u}^{\ub}\frac{1}{|\ub'|^{3+\frac{3}{2}\alpha}} \cdot \delta^{2l}\|\Lb{Z^3}{\phi}\|^2_{L^{2}(\Cb_{\ub'})} d\ub' + \delta^{\frac{1}{4}}\int_{0}^{u}\frac{1}{\delta}\cdot\delta^{2l-1} \||\ub|^{\frac{\alpha}{2}}L{Z^3}{\phi}\|^2_{L^2(C_{u'})} d u'\\
 &\lesssim \delta^{\frac{1}{4}}M^2.
\end{align*}

For $S_{3}$, we have
\begin{align*}
 S_3 &\lesssim   \delta^{2l-1}\doubleint_{\mathcal{D}_{\ub,u}} \big( |\ub|^{-\frac{3}{2}-\frac{\alpha}{4}} \delta^{\frac{1}{4}}M\big)|\nablaslash{Z^3}{\phi}||L{Z^3} {\phi}|\\
&\lesssim M \delta^{-\frac{1}{4}}\int_{0}^{u} \big(\delta^{l} \|\nablaslash{Z^3}{\phi}\|_{L^2(C_{u'})}\big)\big(\delta^{l-\frac{1}{2}} \||\ub|^{\frac{\alpha}{2}}L{Z^3}{\phi}\|_{L^2(C_{u'})}\big) d u'\\
 &\lesssim \delta^{\frac{3}{4}}M^3.
\end{align*}

For $S_{4}$, we have
\begin{align*}
 S_4 &\lesssim   \delta^{2l-1}\doubleint_{\mathcal{D}_{\ub,u}} \big( |\ub|^{-1+\alpha} \delta^{-\frac{1}{2}}I_3\big)|\nablaslash{Z^3}{\phi}||L{Z^3} {\phi}|\\
&\lesssim  \delta^{-1}\int_{0}^{u}  \big(\delta^{l} \|\nablaslash{Z^3}{\phi}\|_{L^2(C_{u'})}\big)\big(\delta^{l-\frac{1}{2}} \||\ub|^{\frac{\alpha}{2}}L{Z^3}{\phi}\|_{L^2(C_{u'})}\big) d u'.
\end{align*}
By virtue of 
\eqref{Eenergy estimates top Lb}, we can bound $\delta^{l} \|\nablaslash{Z^3}{\phi}\|_{L^2(C_{u'})}$ to derive
\begin{align*}
 S_4 &\lesssim  C(I_4) M.
\end{align*}

The estimates on $S_1$, $S_2$, $S_3$ and $S_4$ together yield
\begin{equation}\label{3rd order S Lb}
S \lesssim \delta^{\frac{1}{4}}M^3 + C(I_4) M.
\end{equation}

\bigskip

For $T$, according to \eqref{null form bound}, we have
\begin{align*}
 T &\lesssim \sum_{p+q \leq 3,\atop p\leq 2, q\leq 2} \delta^{2l-1}\doubleint_{D_{\ub,u}} |\ub|^\alpha |\partial Z^{p}\phi| |\partial_g  Z^{q}\phi||L Z^3 \phi|.
\end{align*}
where $\partial\in\{\nablaslash, \Lb\}$ and  $\partial_g \in\{\nablaslash, L\}$.

We first consider $\partial_{g}=\slashed{\nabla}, \partial=\Lb$, and denote its contribution by $T_{1}$, then we see all the other cases are lower order compared to this case. By \eqref{L4 Sobolev} we have:
\begin{align*}
T_{1}\lesssim\delta^{2l-1}\sum_{p+q \leq 3,\atop p\leq 2, q\leq 2} \int_{1-u}^{\ub}\int_{0}^{u}\ub^{\prime\alpha/2}\|\Lb Z^{p}\phi\|_{L^{4}(S_{u',\ub'})}\|\slashed{\nabla}Z^{q}\phi\|_{L^{4}(S_{u',\ub'})}\|\ub^{\prime\alpha/2}LZ^{3}\phi\|_{L^{2}(S_{u',\ub'})}du'd\ub'
\end{align*}
By the second of \eqref{L2 Sobolev},
\begin{align*}
\delta^{l'}\|\Lb Z^{p}\phi\|_{L^{4}(S_{u',\ub'})}&\lesssim\delta^{l'}\ub^{\prime-1/2}\|\Lb Z_{b}Z^{p}\phi\|_{L^{2}(\Cb_{\ub'})}^{1/2}\cdot\\
&\left(\|\Lb Z^{p}\phi\|_{L^{2}(\Cb_{\ub'})}+\|\Lb Z_{g}Z^{p}\phi\|_{L^{2}(\Cb_{\ub'})}\right)^{1/2}\\
&\lesssim\delta^{-1/2}\ub^{\prime-1/2}\left(I_{4}+\delta^{1/8}M^{3/2}\right)\lesssim\delta^{-1/2}
\ub^{\prime-1/2}I_{4}
\end{align*}
provided that $\delta$ is sufficiently small.

On the other hand, by \eqref{L4 Sobolev},
\begin{align*}
\delta^{l''}\|\slashed{\nabla}Z^{q}\phi\|_{L^{4}(S_{\ub',u'})}\lesssim\ub^{\prime-1/2}\left(\delta^{l''}\|\slashed{\nabla}Z^{q}\phi\|_{L^{2}(S_{\ub',u'})}+\delta^{l''}\|\slashed{\nabla}Z_{g}Z^{q}\phi\|_{L^{2}(S_{\ub',u'})}\right)
\end{align*}
Therefore we have:
\begin{align*}
T_{1}&\lesssim\delta^{l-1}\sum_{p+q \leq 3,\atop p\leq 2, q\leq 2}\int_{0}^{u}\delta^{-1/2}I_{4}\left(\delta^{l''}\|\slashed{\nabla}Z^{q}\phi\|_{L^{2}(C_{u'})}+\delta^{l''}\|\slashed{\nabla}Z_{g}Z^{q}\phi\|_{L^{2}(C_{u'})}\right)\\
&\cdot\|\ub^{\prime\alpha/2}LZ^{3}\phi\|_{L^{2}(C_{u'})}du'\lesssim\delta^{l-1}\sum_{p+q \leq 3,\atop p\leq 2, q\leq 2}\int_{0}^{u}\delta^{-1/2}I_{4}^{2}\|\ub^{\prime\alpha/2}LZ^{3}\phi\|_{L^{2}(C_{u'})}du'
\end{align*}
By Cauchy-Schwarz, this implies:
\begin{equation*}
T_{1}\lesssim\delta^{-1}\sum_{p+q \leq 3,\atop p\leq 2, q\leq 2}\int_{0}^{u}\delta^{2l-1}\|\ub^{\prime\alpha/2}LZ^{3}\phi\|^{2}_{L^{2}(C_{u'})}du'+I_{4}^{4}
\end{equation*}
If $\partial=\Lb$, $\partial_{g}=L$, we denote its contribution by $T_{2}$, then by the estimates we have derived for $\|\ub^{\prime\alpha/2}LZ^{q}\phi\|_{L^{2}(C_{u})}$, a similar argument leads to the estimate on $T_{2}$:
\begin{align*}
T_{2}\lesssim\delta^{-1/2}\sum_{p+q \leq 3,\atop p\leq 2, q\leq 2}\int_{0}^{u}\delta^{2l-1}\|\ub^{\prime\alpha/2}LZ^{3}\phi\|^{2}_{L^{2}(C_{u'})}du'+\delta^{1/2}I_{4}^{4}
\end{align*}
If $\partial=\slashed{\nabla}$, the estimates for $\partial_{g}Z^{q}\phi$ are the same as before. While for $\slashed{\nabla}Z^{p}\phi$, we have, if $\delta$ is sufficiently small:
\begin{align*}
\delta^{l''}\|\ub^{\prime\alpha/2}\slashed{\nabla}Z^{p}\phi\|_{L^{4}(S_{u'\ub'})}&\lesssim \delta^{l''}\ub^{\prime-1+\alpha/4}\|\Lb Z^{p}Z_{g}\phi\|_{L^{2}(\Cb_{\ub'})}^{1/2}\\
&\cdot\left(\|\ub^{\prime\alpha/2}\slashed{\nabla}Z^{p}\phi\|_{L^{2}(\Cb_{\ub'})}+\|\ub^{\prime\alpha/2}\slashed{\nabla}Z_{g}Z^{p}\phi\|_{L^{2}(\Cb_{\ub'})}\right)^{1/2}\\
&\lesssim\delta^{1/4}I_{4}
\end{align*}
This bound is better than that of $\delta^{l''}\|\Lb Z^{p}\phi\|_{L^{4}(S_{u'\ub'})}$. Therefore finally we obtain:
\begin{align}\label{3rd order T Lb}
T\lesssim\delta^{-1}\sum_{p+q \leq 3,\atop p\leq 2, q\leq 2}\int_{0}^{u}\delta^{2l-1}\|\ub^{\prime\alpha/2}LZ^{3}\phi\|^{2}_{L^{2}(C_{u'})}du'+I_{4}^{4}
\end{align}

\bigskip

It remains to control $W = \delta^{2l-1}\int\!\!\!\int_{\mathcal{D}_{\ub,u}} \frac{|\ub|^{\alpha}}{r}|L{Z^3} {\phi}||\Lb{Z^3}{\phi}|$. By Cauchy-Schwarz inequality, we have
\begin{align*}
W &\lesssim \delta^{2l-1}\Big( \doubleint_{\mathcal{D}_{\ub,u}}\frac{\delta}{|\ub'|^{2-\alpha}}|\Lb{Z^3}{\phi}|^2+  \doubleint_{\mathcal{D}_{\ub,u}} \frac{1}{\delta}||\ub'|^{\frac{\alpha}{2}}L{Z^3}{\phi}|^2\Big)\\
&= \delta^{2l}\int_{1-u}^{\ub}\frac{1}{|\ub'|^{2-\alpha}} \|\Lb{Z^3}{\phi}\|^2_{L^{2}(\Cb_{\ub'})} d\ub' + \delta^{2l-2}\int_{0}^{u}\||\ub|^{\frac{\alpha}{2}}L{Z^3}{\phi}\|^2_{L^2(C_{u'})} d u'.
\end{align*}
We use \eqref{Eenergy estimates top Lb} to bound the first term in the last line. Since $\alpha <1$, for sufficiently small $\delta$, we obtain
\begin{equation}\label{3rd order W Lb}
W \lesssim I_{4}^2 +  \delta^{-1} \int_{0}^{u} \delta^{2l-1} \||\ub|^{\frac{\alpha}{2}}L{Z^3}{\phi}\|^2_{L^2(C_{u'})} d u'.
\end{equation}

By combining \eqref{3rd order S Lb}, \eqref{3rd order T Lb} and \eqref{3rd order W Lb}, we obtain
\begin{align*}
& \ \ \delta^{2l-1}\int_{C_{u}}|\ub|^{\alpha} |L {Z^3}{\phi}|^{2}+\delta^{2l-1} \int_{\underline{C}_{\underline{u}}}|\ub|^{\alpha}|\nablaslash {Z^3} {\phi}|^{2}\\
&\lesssim I_4^4 + \delta^{\frac{1}{4}}M^3 + C(I_3)M+\delta^{-1} \int_{0}^{u} \delta^{2l-1} \||\ub|^{\frac{\alpha}{2}}L{Z^3}{\phi}\|^2_{L^2(C_{u'})} d u'
\end{align*}
The last term on the right-hand side can be removed by the Gronwall's inequality so that
\begin{align*}
& \ \ \delta^{2l-1}\int_{C_{u}}|\ub|^{\alpha} |L {Z^3}{\phi}|^{2}+\delta^{2l-1} \int_{\underline{C}_{\underline{u}}}|\ub|^{\alpha}|\nablaslash {Z^3} {\phi}|^{2}\\
&\lesssim I_4^4 + \delta^{\frac{1}{4}}M^3 + C(I_3)M.
\end{align*}

This finally proves that, for all $0\leq l \leq 3$, we have
\begin{equation}\label{top Eenergy estimates 3 L}
\delta^{l-\frac{1}{2}}\||\ub|^\frac{\alpha}{2}L Z_b^l Z_g^{3-l} \phi\|_{L^2(C_u)} +\delta^{l-\frac{1}{2}}\|\nablaslash Z_b^l Z_g^{3-l} \phi\|_{L^2(\Cb_{\ub})} \lesssim C(I_4)+ C(I_3)M^\frac{1}{2}+ \delta^{\frac{1}{8}}M^{\frac{3}{2}}.
\end{equation}

By combining this estimates with \eqref{top Eenergy estimates 3 L} and \eqref{energy estimates k leq 2}, for sufficiently small $\delta$, this finally proves
\begin{equation}\label{energy estimates k leq 3}
E_{\leq 3}(u,\ub) + \Eb_{\leq 3}( u,\ub) \leq C(I_4).	
\end{equation}
This is the end of the bootstrap argument and the Proposition \ref{proposition bootstrap} has been proved.

\subsection{Higher Order Estimates}
This subsection is devoted to prove a higher order analogue of Proposition \ref{proposition bootstrap}:
\begin{proposition}\label{proposition higher order estimates}
Given a positive integer $n \geq 3$, there exists $\delta_0 >0$, so that for all $\delta < \delta_0$,  for all $u \in (0,u^*)$ and $\ub \in (1-u^*, \ub^*)$, we have
\begin{equation}\label{higher order estimates}
E_{\leq n}(u,\ub) + \Eb_{\leq n}( u,\ub) \leq C(I_{n+1}),
\end{equation}
and
\begin{equation}\label{higher order L infty esimates}
\begin{split}
 \|\nablaslash Z_b^l Z_g^{k-l}\|_{L^\infty(S_{\ub,u})} &\lesssim  \delta^{\frac{1}{4}-l}|\ub|^{-\frac{3}{2}-\frac{\alpha}{4}}C(I_{n+1}), \ \ 0 \leq l \leq k \leq n-3,\\
 \|L Z_b^l Z_g^{k-l} \phi\|_{L^{\infty}(S_{\ub,u})}&\lesssim \delta^{\frac{1}{2}-l}|\ub|^{-2}C(I_{n+1}),  \ \ 0 \leq l \leq k \leq n-3,\\
\|\Lb Z_b^l Z_g^{k-l} \phi\|_{L^\infty(S_{\ub,u})} &\lesssim \delta^{-\frac{1}{2}-l}{|\ub|}^{-1} C(I_{n+1}), \ \ 0 \leq l \leq k \leq n-3.
\end{split}
\end{equation}
where $C(I_{n+1})$ is a constant depending only on $I_{n+1}$.
\end{proposition}
\begin{remark}
Although $C(I_n)$ and $\delta_0$ in the proposition may depend on the integer $n$, in the rest of the paper, we only need the result for $n = 12$.
\end{remark}

We prove \eqref{higher order estimates} and \eqref{higher order L infty esimates} together by induction on $n$. For $n=3$, the proposition has been achieved in the previous subsection. For $n\geq 4$, we assume that the proposition holds for all $n'$ so that $n'\leq n-1$. To prove for $n$, we first make the following bootstrap assumption: We choose a large constant $M$, so that
\begin{equation}\label{higher order bootstrap assumption L 2}
E_{n}(u,\ub) + \Eb_{n}( u,\ub) \lesssim M,
\end{equation}
for all $u \in (0,u^*)$ and $\ub \in (1-u^*, \ub^*)$. We remark that $M$ may depend on $\phi$ at the moment. We will show that, if $\delta$ is sufficiently small, then we can make $M$ depend only on $I_{n+1}$. We also remark that the induction hypothesis is
\begin{equation}\label{higher order induction hypothesis L 2}
E_{\leq n-1}(u,\ub) + \Eb_{\leq n-1}( u,\ub) \lesssim C(I_{n}),
\end{equation}
and
\begin{equation}\label{higher order induction hypothesis L infty}
\begin{split}
 \|\nablaslash Z_b^l Z_g^{k-l}\phi\|_{L^\infty(S_{\ub,u})} &\lesssim  \delta^{\frac{1}{4}-l}|\ub|^{-\frac{3}{2}-\frac{\alpha}{4}}C(I_n), \ \ 0 \leq l \leq k \leq n-4,\\
 \|L Z_b^l Z_g^{k-l} \phi\|_{L^{\infty}(S_{\ub,u})}&\lesssim \delta^{\frac{1}{2}-l}|\ub|^{-2}C(I_n),  \ \ 0 \leq l \leq k \leq n-4,\\
\|\Lb Z_b^l Z_g^{k-l} \phi\|_{L^\infty(S_{\ub,u})} &\lesssim \delta^{-\frac{1}{2}-l}{|\ub|}^{-1} C(I_n), \ \ 0 \leq l \leq k \leq n-4.
\end{split}
\end{equation}
for all $u \in (0,u^*)$ and $\ub \in (1-u^*, \ub^*)$.

We claim that, together with the induction hypothesis \eqref{higher order induction hypothesis L 2} and \eqref{higher order induction hypothesis L infty}, the bootstrap assumption \eqref{higher order bootstrap assumption L 2} implies
\begin{equation}\label{higher order bootstrap assumption L infty}
\begin{split}
 \|\nablaslash Z_b^l Z_g^{n-2-l}\|_{L^\infty(S_{\ub,u})} &\lesssim  \delta^{\frac{1}{4}-l}|\ub|^{-\frac{3}{2}-\frac{\alpha}{4}}M, \ \ 0 \leq l \leq n-3,\\
 \|L Z_b^l Z_g^{n-2-l} \phi\|_{L^{\infty}(S_{\ub,u})}&\lesssim \delta^{\frac{1}{2}-l}|\ub|^{-2}M,  \ \ 0 \leq l \leq  n-3,\\
\|\Lb Z_b^l Z_g^{n-2-l} \phi\|_{L^\infty(S_{\ub,u})} &\lesssim \delta^{-\frac{1}{2}-l}{|\ub|}^{-1} M, \ \ 0 \leq l \leq n-3.
\end{split}
\end{equation}

The bound on $\|\nablaslash Z_b^l Z_g^{n-3-l}\|_{L^\infty(S_{\ub,u})}$ is straightforward: we simply use Sobolev inequalities by affording two more $\Omega_{ij}$ derivatives. The derivation is exactly the same as for \eqref{p_7}.

The bound on $\|L Z_b^l Z_g^{n-3-l} \phi\|_{L^{\infty}(S_{\ub,u})}$ relies on the \eqref{commuted null frame equation}, i.e.
\begin{equation}\label{E 8}
\Box Z^{n-3} \phi = \sum_{p+q \leq n-3}Q(\nabla Z^p \phi, \nabla Z^q \phi).
\end{equation}
According to the structure of null forms, we can rewrite it as the following inequality:
\begin{equation}\label{E 9}
\Lb|L Z^{n-3}\phi|\lesssim a|L Z^{n-3} \phi|+b,
\end{equation}
where
\begin{equation*}
 a=\frac{1}{r}+
 \big(|\Lb \phi|+|\slashed{\nabla} \phi|\big),
\end{equation*}
 and
\begin{equation*}
 b=|\slashed{\triangle} Z^{n-3} \phi|+\frac{1}{r}|\Lb Z^{n-3} \phi|+\sum_{p+q \leq n-3}\big(|\Lb Z^p \phi||\slashed{\nabla} Z^q\phi|+|\slashed{\nabla} Z^p\phi||\slashed{\nabla} Z^q\phi| \big).
\end{equation*}
We claim that
\begin{equation}\label{bound on a and b}
\begin{split}
 \|a\|_{L^1_u L^\infty(S_{\ub,u})} &\lesssim |\ub|^{-1}\delta^{-\frac{1}{2}-l} M,\\
  \|b\|_{L^1_u L^\infty(S_{\ub,u})} &\lesssim |\ub|^{-2}\delta^{-\frac{1}{2}-l} M.
\end{split}
\end{equation}
To prove this claim, we first notice that all the terms have already been bounded by the induction hypothesis except for the top order terms, i.e. $|\slashed{\triangle} Z^{n-3} \phi|$, $\frac{1}{r}|\Lb Z^{n-3} \phi|$, $|\Lb Z^{n-3} \phi||\slashed{\nabla} \phi|$, $|\Lb \phi||\slashed{\nabla} Z^{n-3}\phi|$ and $|\slashed{\nabla}\phi||\slashed{\nabla} Z^{n-3}\phi|$ appeared in $b$. In view of the bound on 
$\|\nablaslash Z_b^l Z_g^{n-3-l}\phi\|_{L^\infty(S_{\ub,u})}$ derived above, it suffices to bound $|\Lb Z^{n-3} \phi|$. According to Sobolev inequality, we have
\begin{align*}
 \|\Lb Z^{n-3} \phi\|_{L^1_u L^\infty(S_{\ub,u})}&\lesssim |\ub|^{-1}\sum_{0\leq j \leq 2} \|\Omega^j \Lb  Z^{n-3}  \phi\|_{L^1_u L^2(S_{\ub,u})}\\
 &\lesssim |\ub|^{-1}\delta^{\frac{1}{2}}\sum_{0\leq j \leq 2} \|\Lb \Omega^j  Z^{n-3} \phi\|_{L^2_u L^2(S_{\ub,u})}\\
 &\lesssim |\ub|^{-1}\delta^{\frac{1}{2}} M.
\end{align*}
We thus proved \eqref{bound on a and b}. By virtue of Gronwall's inequality, \eqref{E 9} yields the desired estimates for $L Z_b^l Z_g^{k-l} \phi$ in \eqref{higher order bootstrap assumption L infty}.

The estimates on $\|\Lb Z_b^l Z_g^{n-3-l} \phi\|_{L^\infty(S_{\ub,u})}$ relies on the use of equation \eqref{commuated main equatin for energy estimates}. In fact, we have
\begin{equation*}
\Box Z^{n-3} \phi = \sum_{p+q \leq n-3}Q(\nabla Z_g^p \phi, \nabla Z_g^q \phi).
\end{equation*}
Let $y = |\ub||\Lb Z_b^l Z_g^{n-3-l} \phi|$. By computing $L\big(\ub^{2}(\Lb Z_b^l Z_g^{n-3-l} \phi)^{2}\big)$, we have
\begin{equation*}
L y^{2}\lesssim \big(\frac{\delta}{|\ub|^{2}}+ \frac{\delta^\frac{1}{4}}{|\ub|^{2}}M\big) y^{2}+ \frac{\delta^{\frac{1}{4}-l}}{|\ub|^{2}} M y.
\end{equation*}
By integrating this equation, we obtain
\begin{equation}\label{b1}
\big||\ub||\Lb Z_b^l Z_g^{n-3-l}  \phi|(\ub,u,\theta)-C|1-u||\Lb Z_b^l Z_g^{n-3-l} \phi|(1-u,u,\theta)\big| \lesssim  \delta^{-\frac{1}{4}-l}M,
\end{equation}
Therefore, according to \eqref{L2 Sobolev initial}, we finally obtain
\begin{equation*}
\|\Lb Z_b^l Z_g^{n-3-l}\phi\|_{L^\infty(S_{\ub,u})} \lesssim \frac{\delta^{-\frac{1}{2}-l}}{|\ub|} I_{n+1} + \delta^{-\frac{1}{4}-l}|\ub|^{-1}M.
\end{equation*}

\bigskip

To finish the proof of Proposition \ref{higher order estimates}, it remains to improve the constant $M$ in \eqref{higher order bootstrap assumption L 2}. The procedure is exactly the same as for the proof of $E_{3}(u,\ub)$ and $\Eb_3(u,\ub)$ in previous subsection.

\bigskip

We 
replace $\phi$ by $Z^n \phi$ and take $X = \Lb$ in \eqref{fundamental energy identity} and \eqref{commuated main equatin for energy estimates}, this yields
\begin{equation}\label{E n Lb}
\begin{split}
\delta^{2l}\int_{C_{u}}|\nablaslash Z^n \phi|^{2}&+\delta^{2l}\big|\int_{\underline{C}_{\underline{u}}}|\Lb Z^n \phi|^{2} \lesssim I_{n+1}^2 + \delta^{2l}\big|\doubleint_{D_{\ub,u}}Q(\nabla Z^n \phi, \nabla \phi)\Lb Z^n \phi \big| \\
&  + \sum_{p+q \leq n,\atop p<n, q<n} \delta^{2l}\big|\doubleint_{D_{\ub,u}}
Q(\nabla Z^{p}\phi, \nabla Z^{q}\phi)\Lb Z^n \phi \big|+\delta^{2l}\big|\doubleint_{D_{\ub,u}} \frac{1}{r} \Lb Z^n \phi \cdot L Z^n\phi\big|.
\end{split}
\end{equation}
We rewrite the right-hand side $I_{n+1}^2 +  S + T +W$, where $S$, $T$ and $W$ denote the three bulk integral terms in \eqref{E n Lb}. We will bound $S$, $T$ and $W$ one by one.

\bigskip

We start with $S$ which bounded by the sum of the following terms:
\begin{align*}
S_1 &= \delta^{2l}\doubleint_{\mathcal{D}_{\ub,u}}|\Lb {\phi}|\big(|L Z^n {\phi}|+|\nablaslash Z^n{\phi}|\big)|\Lb Z^n {\phi}|, \\
S_2 &= \delta^{2l}\doubleint_{\mathcal{D}_{\ub,u}}|L{\phi}|\big(|\nablaslash Z^n {\phi}|+|\Lb Z^n {\phi}|\big)|\Lb Z^n {\phi}|,\\
S_3 &= \delta^{2l}\doubleint_{\mathcal{D}_{\ub,u}}|\nablaslash {\phi}|\big(|L Z^n {\phi}|+|\Lb Z^n {\phi}|+|\nablaslash Z^n {\phi}|\big)|\Lb Z^n{\phi}|.
\end{align*}

In view of the forms of $S_1$, $S_2$ and $S_3$ appeared in the subsection for the estimates on $E_{3}(u,\ub)$ and $\Eb_3(u,\ub)$, i.e. the derivation of the inequalities \eqref{3rd order S1 L}, \eqref{3rd order S2 L} and \eqref{3rd order S3 L}, we can proceed \emph{exactly} in the same way (simply replace all the $Z^3\phi$ by $Z^n \phi$). We take $S_1$ as an example to illustrate the process: by the bound on $\Lb{\phi}$ in $L^\infty(S_{\ub,u})$, we have
\begin{align*}
S_{1}&\lesssim \delta^{2l-\frac{1}{2}}\doubleint_{\mathcal{D}_{\ub,u}} |\ub|^{-1}|\nablaslash Z^n{\phi}| |\Lb Z^n{\phi}|+ \delta^{2l-\frac{1}{2}}\doubleint_{\mathcal{D}_{\ub,u}} |\ub|^{-\frac{1+\alpha}{2}} \cdot \big(|\ub|^{\frac{\alpha}{2}} L Z^n{\phi}|\big)\cdot\big(|\ub|^{-\frac{1}{2}}|\Lb Z^n{\phi}|\big)\\
&=S_{11}+S_{12}.
\end{align*}
We bound $S_{11}$ exactly as the derivation for \eqref{S11}:
\begin{align*}
S_{11}& \lesssim \delta^{2l-1}\doubleint_{\mathcal{D}_{\ub,u}}|\nablaslash Z^n{\phi}|^{2}+\delta^{2l}\doubleint_{\mathcal{D}_{\ub,u}}|\ub|^{-2}|\Lb\Omega^{3}{\phi}|^{2}\\
&\lesssim \int_{0}^{u}\frac{1}{\delta} \cdot \delta^{2l}\int_{C_{u'}}|\nablaslash Z^n{\phi}|^{2} du'+\int_{1-u}^{\ub} |\ub'|^{-2} \delta^{2l}\int_{\Cb_{\ub}}|\Lb Z^n{\phi}|^{2} d\ub'.
\end{align*}
Similarly, we have
\begin{align*}
S_{12} \lesssim \delta^{\frac{1}{4}}M^2.
\end{align*}

We give the final result on as follows:
\begin{equation}\label{top order S Lb}
S \lesssim \delta^{\frac{1}{4}}M^2 + \int_{0}^{u}\frac{1}{\delta} \cdot \delta^{2l}\int_{C_{u'}}|\nablaslash Z^n{\phi}|^{2} du'+\int_{1-u}^{\ub} |\ub'|^{-2} \delta^{2l}\int_{\Cb_{\ub}}|\Lb Z^n{\phi}|^{2} d\ub'.
\end{equation}

\bigskip

For $T$, we have
\begin{align*}
 T &\lesssim \sum_{p+q \leq n,\atop p\leq n-1, q\leq n-1} \delta^{2l}\doubleint_{D_{\ub,u}} |\partial Z^{p}\phi| |\partial_g  Z^{q}\phi||\Lb Z^n \phi|.
\end{align*}
where $\partial\in\{\nablaslash, \Lb\}$ and  $\partial_g \in\{\nablaslash, L\}$. The estimate for $T$ follows exactly the same as we derive \eqref{T Lb}.
we have:
\begin{align}\label{top order T Lb}
T\lesssim\delta^{1/4}M^{3}
\end{align}
\bigskip

For $W = \delta^{2l}\int\!\!\!\int_{\mathcal{D}_{\ub,u}} \frac{1}{r}|L{Z^n} {\phi}||\Lb{Z^n}{\phi}|$, we have
\begin{align*}
W &{\lesssim} \int_{1-u}^{\ub}\frac{1}{|\ub'|^{2-\alpha}} \cdot \delta^{2l+1}\|\Lb{Z^n}{\phi}\|^2_{L^{2}(\Cb_{\ub'})} d\ub' + \delta^{2l-1}\int_{0}^{u}\||\ub|^{\frac{\alpha}{2}}L{Z^n}{\phi}\|^2_{L^2(C_{u'})} d u'\\
&\lesssim \delta M^2.
\end{align*}

\bigskip

The estimates on $S$, $T$ and $W$, together with \eqref{E n Lb}, imply that
\begin{align*}
\delta^{2l}\int_{C_{u}}|\nablaslash Z^n \phi|^{2}+\delta^{2l}\int_{\underline{C}_{\underline{u}}}|\Lb Z^n \phi|&^{2} \lesssim I_{n+1}^2 + \delta^{\frac{1}{4}}M^3 \\
&\ \ +\int_{0}^{u}\frac{1}{\delta} \cdot \delta^{2l}\int_{C_{u'}}|\nablaslash Z^n{\phi}|^{2} du'+\int_{1-u}^{\ub} |\ub'|^{-2} \delta^{2l}\int_{\Cb_{\ub}}|\Lb Z^n{\phi}|^{2} d\ub'.
\end{align*}
The last two terms can be removed by Gronwall's inequality. Therefore, we obtain
\begin{equation}\label{top order Energy Lb}
\delta^{l}\|\nablaslash Z_b^l Z_g^{n-l} \phi\|_{L^2(C_u)} +\delta^{l}\|\Lb Z_b^l Z_g^{n-l} \phi\|_{L^2(\Cb_{\ub})} \lesssim I_{n+1} + \delta^{\frac{1}{8}}M^{\frac{3}{2}}.
\end{equation}

\bigskip

We now change the multiplier vector field to $u^\alpha L$ to derive
\begin{equation*}
\begin{split}
\delta^{2l-1}\int_{C_{u}}|\ub|^{\alpha}&|L {Z^n}{\phi}|^{2}+\delta^{2l-1} \int_{\underline{C}_{\underline{u}}}|\ub|^{\alpha}|\nablaslash {Z^n} {\phi}|^{2}\lesssim I_{n+1}^2  +\delta^{2l-1}\big|\doubleint_{\mathcal{D}_{\ub,u}}|\ub|^{\alpha} Q(\nabla Z^n \phi, \nabla\phi)L {Z^n} {\phi}\big| \\
& + \sum_{p+q \leq n,\atop p\leq n-1, q \leq n-1} \delta^{2l-1} \big|\doubleint_{D_{\ub,u}} |\ub|^{\alpha} Q(\nabla Z^{p}\phi, \nabla Z^{q}\phi)L {Z^n}{\phi}\big| +\delta^{2l-2}\big|\doubleint_{D_{\ub,u}} \frac{|\ub|^{\alpha}}{r} \Lb {Z^n} {\phi} \cdot L {Z^n}{\phi}\big|,
\end{split}
\end{equation*}
We rewrite the above inequality as
\begin{equation}\label{E top LL}
\delta^{2l-1}\int_{C_{u}}|\ub|^{\alpha}|L {Z^n}{\phi}|^{2}+\delta^{2l-1} \int_{\underline{C}_{\underline{u}}}|\ub|^{\alpha}|\nablaslash {Z^n} {\phi}|^{2}\lesssim I_{n+1}^2  +  S + T +W.
\end{equation}
where $S$, $T$ and $W$ denote the three bulk integral terms. We bound $S$, $T$ and $W$ one by one.

To bound $S$, we can follow exactly the same way as the derivation for \eqref{3rd order S Lb} (we simply replace all the $Z^3\phi$'s by $Z^n \phi$), this gives
\begin{equation*}
S \lesssim \delta^{\frac{1}{4}}M^3 + C(I_{n+1}) M.
\end{equation*}

\bigskip

To bound $T$, we can follow exactly the same way as the derivation for \eqref{3rd order T Lb} (we simply replace all the $Z^3\phi$'s by $Z^n \phi$ and $Z^2 \phi$'s by $Z^{n-1}\phi$). 
We the obtain
\begin{equation*}
T\lesssim
\delta^{-1}\sum_{p+q\leq n,\atop p\leq n-1,q\leq n-1}\int_{0}^{u}\delta^{2l-1}\|\ub^{\prime\alpha/2}LZ^{n}\phi\|^{2}_{L^{2}(C_{u'})}du'+I_{n+1}^{4}
\end{equation*}

\bigskip

To bound $W$, we can follow exactly the same way as the derivation for \eqref{3rd order W Lb} by replacing the $Z^3\phi$'s by $Z^n \phi$, this gives
\begin{equation*}
W \lesssim I_{n+1}^2 +  \delta^{-1} \int_{0}^{u} \delta^{2l-1} \||\ub|^{\frac{\alpha}{2}}L{Z^n}{\phi}\|^2_{L^2(C_{u'})} d u'.
\end{equation*}

\bigskip
The estimates on $S$, $T$ and $W$, together with \eqref{E n Lb}, imply that
\begin{align*}
& \ \ \delta^{2l-1}\int_{C_{u}}|\ub|^{\alpha} |L {Z^n}{\phi}|^{2}+\delta^{2l-1} \int_{\underline{C}_{\underline{u}}}|\ub|^{\alpha}|\nablaslash {Z^n} {\phi}|^{2}\\
&\lesssim I_{n+1}^4 + \delta^{\frac{1}{4}}M^3 + C(I_n)M+\delta^{-1} \int_{0}^{u} \delta^{2l-1} \||\ub|^{\frac{\alpha}{2}}L{Z^n}{\phi}\|^2_{L^2(C_{u'})} d u'
\end{align*}
By the Gronwall's inequality again, for $l\leq n$, we finally obtain
\begin{equation}\label{top order Energy L}
\delta^{l-\frac{1}{2}}\||\ub|^\frac{\alpha}{2}L Z_b^l Z_g^{n-l} \phi\|_{L^2(C_u)} +\delta^{l-\frac{1}{2}}\|\nablaslash Z_b^l Z_g^{n-l} \phi\|_{L^2(\Cb_{\ub})} \lesssim C(I_{n+1})+ C(I_n)M^\frac{1}{2}+ \delta^{\frac{1}{8}}M^{\frac{3}{2}}.
\end{equation}

For sufficiently small $\delta$, the estimate \eqref{top order Energy Lb} and \eqref{top order Energy L} show that
\begin{equation*}\label{energy estimates k leq n}
E_{n}(u,\ub) + \Eb_{n}( u,\ub) \leq C(I_{n+1}).	
\end{equation*}
This completes the bootstrap argument and the Proposition \ref{proposition higher order estimates} has been proved.

\subsection{Existence based on \textit{a priori} estimates}

The existence of solutions of \eqref{Main Equation} follows immediately from the a priori energy estimates derived previously. Since the procedure is standard, we only give a sketch of the proof in this subsection.

We start with solving local solution for Cauchy problem with data prescribed on $\Sigma_{1}$ with $1-\delta \leq r \leq 1$. Therefore, we obtain a local solution confined in the region bounded by $C_\delta$ and $\Cb_1$. In particular, on a neighborhood of $S_{1,0}$ on the incoming cone $\Cb_1$, the solution has been constructed.

We then use $C_0$ and $\Cb_1$ as initial hypersurfaces. The classical local existence result \cite{R-90} of Rendall can be applied in this situation. Therefore, we know that there exists a solution in the entire spacetime neighborhood (which lies in the domain of dependence of $C_0$ and $\Cb_1$) of $S_{1, 0}$.Combined with the local solution of the Cauchy problem, we have constructed a local solution for $t\in [1, 1+\epsilon]$ for some small $\epsilon$.

Since the a priori energy estimates (as well as the companying $L^\infty$ estimates) depends only on the size $I_{n}$ of the rescaled data on $\Sigma_1$,  this solution is well behaved on $\Sigma_{1+\epsilon}$. Therefore, we can use this as initial surface (instead of $\Sigma_1$) to repeat the above argument. Eventually, we obtain an solution in the entire short pulse region.

\subsection{Improved Estimates on $C_\delta$}
Recall that, given $n \geq 12$, in the short pulse region, we have derived the following a priori $L^{\infty}$ estimates on the solution $\phi$:
\begin{equation*}
\begin{split}
 \|\nablaslash Z_b^l Z_g^{k-l} \phi\|_{L^\infty(S_{\ub,u})} &\lesssim  \delta^{\frac{1}{4}-l}|\ub|^{-\frac{3}{2}-\frac{\alpha}{4}}C(I_{n+1}), \ \ 0 \leq l \leq k \leq n-3,\\
 \|L Z_b^l Z_g^{k-l} \phi\|_{L^{\infty}(S_{\ub,u})}&\lesssim \delta^{\frac{1}{2}-l}|\ub|^{-2}C(I_{n+1}),  \ \ 0 \leq l \leq k \leq n-3,\\
\|\Lb Z_b^l Z_g^{k-l} \phi\|_{L^\infty(S_{\ub,u})} &\lesssim \delta^{-\frac{1}{2}-l}{|\ub|}^{-1} C(I_{n+1}), \ \ 0 \leq l \leq k \leq n-3.
\end{split}
\end{equation*}
The goal of this section is to improve these bounds for the solution on $C_\delta$. More precisely, we will prove that
\begin{proposition}\label{improved estimates on C delta}
On 
$C_\delta$, for sufficiently small $\delta$, we have
\begin{equation}\label{impoved L infty estmiates}
\begin{split}
 \|\nablaslash Z_b^l Z_g^{k-l}\phi \|_{L^\infty(S_{\ub,\delta})} &\lesssim  \delta^{\frac{1}{4}}|\ub|^{-\frac{3}{2}-\frac{\alpha}{4}}C(I_{n+1}), \ \ 0 \leq l \leq k \leq n-3,\\
 \|L Z_b^l Z_g^{k-l} \phi\|_{L^{\infty}(S_{\ub,\delta})}&\lesssim \delta^{\frac{1}{4}}|\ub|^{-2}C(I_{n+1}),  \ \ 0 \leq l \leq k \leq n-3,\\
\|\Lb Z_b^l Z_g^{k-l} \phi\|_{L^\infty(S_{\ub,\delta})} &\lesssim \delta^{\frac{1}{4}}{|\ub|}^{-1} C(I_{n+1}), \ \ 0 \leq l \leq k \leq n-3.
\end{split}
\end{equation}
\end{proposition}
Notice the power of $\delta$ for $L Z_b^l Z_g^{k-l} \phi$ has been modified to $\frac{1}{2}$. Since $C_0$ is also the outer boundary of the small data region (region I), the smallness (in terms of $\delta$) of the solution stated in the proposition is indispensable for the construction of a global solution in the small data region. As we mentioned in the introduction, the proof relies on the following observation: on the $S_{1-\delta,\delta}$ or equivalently the lower boundary of $C_\delta$, the data are \emph{identically zero}. This is because that the data are compactly supported on $\Sigma_1$ between $S_{1-\delta,\delta}$ and $S_{1,0}$. Therefore, even for bad derivatives of $\phi$, it is small at least initially. The idea of the proof is to integrate along the $L$ direction to show that the smallness indeed propagates.

\bigskip

We use the induction argument on the pair $(l,k)$ ($0 \leq k \leq n-3, 0\leq l\leq k$) to prove \eqref{impoved L infty estmiates}. First of all, we give an order on the set of such pairs: we say that $(l',k') < (l,k)$ if one of the following holds: (1) $k'<k$ or (2) $l'<l$, $k=k'$. We do the induction with respect to this order.

\bigskip

For $(l,k)=(0,0)$, the bounds on $L \phi$ and $\nablaslash \phi$ are clear. It remains to prove that
\begin{equation*}
 \|\Lb \phi\|_{L^\infty(S_{\ub,\delta})} \lesssim \delta^{\frac{1}{4}}{|\ub|}^{-1} C(I_n).
\end{equation*}
Recall that, in \eqref{a1}, we have obtained
\begin{equation*}
\big||\ub||\Lb\phi|(\ub,u,\theta)-C|\ub||\Lb\phi|(1-u,u,\theta)\big| \lesssim  \delta^\frac{1}{4}M.
\end{equation*}
In view of the higher order energy estimates derived in the previous subsection, the constant $M$ should be replaced by $C(I_n)$. Let $u=\delta$, then the second term vanishes on the initial sphere $S_{1-\delta,\delta}$. This gives the desired estimates on $\|\Lb \phi\|_{L^\infty(S_{\ub,\delta})}$.

\bigskip

For $(l,k)=(0,k)$,  we can use \eqref{b1} to obtain the desired estimates in a similar way.

\bigskip

We assume that for all $(l',k') < (l,k)$, we have
\begin{equation*}
\begin{split}
 \|\nablaslash Z_b^{l'} Z_g^{k'-l'}\phi \|_{L^\infty(S_{\ub,\delta})} &\lesssim  \delta^{\frac{1}{4}}|\ub|^{-\frac{3}{2}-\frac{\alpha}{4}}C(I_n), \\
 \|L Z_b^{l'} Z_g^{k'-l'}\phi\|_{L^{\infty}(S_{\ub,\delta})}&\lesssim \delta^{\frac{1}{4}}|\ub|^{-2}C(I_n), \\
\|\Lb Z_b^{l'} Z_g^{k'-l'}\phi\|_{L^\infty(S_{\ub,\delta})} &\lesssim \delta^{\frac{1}{4}}{|\ub|}^{-1} C(I_n).
\end{split}
\end{equation*}
For $(l,k)$, we now reduce the estimates to the above induction hypothesis. Because we have already proved the case for $(l,k)=(0,k)$, so we can assume in addition that $l\geq 1$.

\bigskip

We first bound $\nablaslash Z_b^{l} Z_g^{k-l}\phi$. In fact, we have
\begin{align*}
 |\nablaslash Z_b^{l} Z_g^{k-l}\phi| &\lesssim  \frac{1}{|\ub|}|\Omega Z_b^{l} Z_g^{k-l}\phi|\\
 &\leq \frac{1}{|\ub|}\Big(| Z_b \Omega Z_b^{l-1} Z_g^{k-l}\phi|+|[\Omega,Z_b] Z_b^{l-1} Z_g^{k-l}\phi|\Big).
\end{align*}
Since $\Omega \in \mathcal{Z}_g$, we can use induction hypothesis (since we can reduce $l$), the first term is bounded by
\begin{align*}
 \sum_{\partial \in \{L,\Lb, \slashed{\nabla}\}}\frac{1}{|\ub|}| \partial Z_b^{l-1} Z_g^{k-l+1}\phi| \lesssim \delta^{\frac{1}{4}}|\ub|^{-2}C(I_{n+1}).
\end{align*}
For the second term, notice that $[\Omega,Z_b] = Z_b$, therefore, we have decreased the number $k$ by $1$. According to the induction hypothesis, it is bounded by
\begin{align*}
 \sum_{\partial \in \{L,\Lb, \slashed{\nabla}\}}\frac{1}{|\ub|}| \partial Z_b^{l-1} Z_g^{k-l}\phi| \lesssim \delta^{\frac{1}{4}}|\ub|^{-2}C(I_{n+1}).
\end{align*}
This gives the desired estimates for $\nablaslash Z_b^{l} Z_g^{k-l}\phi$.

\bigskip

We turn to the bound on $L Z_b^{l} Z_g^{k-l}\phi$. Evidently, it is bounded by $\sum_{\partial \in \{L,\Lb, \slashed{\nabla}\}}|L \partial Z_b^{l-1} Z_g^{k-l}\phi|$. When $\partial = \nablaslash$ in the sum, it can be bounded directly by the bound on 
$\nablaslash Z_b^{l-1} Z_g^{k-l+1}\phi$ derived above. Therefore, it suffices to bound $L L Z_b^{l-1} Z_g^{k-l}\phi$ and $L \Lb Z_b^{l-1} Z_g^{k-l}\phi$.

For $L \Lb Z_b^{l-1} Z_g^{k-l}\phi$, according to \eqref{commuated main equatin for energy estimates} (where we use $Z_b^{l-1} Z_g^{k-l}$ as commutator vector field), we have
\begin{align*}
 -L\Lb Z_b^{l-1} Z_g^{k-l} \phi+ \slashed{\triangle} Z_b^{l-1} Z_g^{k-l}\phi =\frac{1}{r}(L Z_b^{l-1} Z_g^{k-l}\phi-\Lb Z_b^{l-1} Z_g^{k-l}\phi) + \sum_{p+q\leq k-1} Q(\nabla Z^p\phi, \nabla Z^q \phi).
\end{align*}
The second term on the left-hand side can be bounded by $\nablaslash Z_b^{l} Z_g^{k-l}\phi$. The terms on the right-hand side are all of lower degrees ($<k$) so that they are bounded by the induction hypothesis. Therefore, we have
\begin{equation*}
 \|L \Lb Z_b^{l-1} Z_g^{k-l}\phi\|_{L^{\infty}(S_{\ub,\delta})} \lesssim \delta^{\frac{1}{4}}|\ub|^{-2}C(I_{n+1}),
\end{equation*}

For $L L Z_b^{l-1} Z_g^{k-l}\phi$, we use the following identity:
\begin{align*}
 L L Z_b^{l-1} Z_g^{k-l}\phi &= L\Big(\frac{1}{\ub}\big(S Z_b^{l-1} Z_g^{k-l}\phi - u\Lb Z_b^{l-1} Z_g^{k-l}\phi\big)\Big).
\end{align*}
The second term on the right-hand side can be bounded directly by the bound on $L \Lb Z_b^{l-1} Z_g^{k-l}\phi$ just derived. Therefore, it suffices to control the contribution from the first term, i.e.
\begin{align*}
 L\Big(\frac{1}{\ub}\big(S Z_b^{l-1} Z_g^{k-l}\phi\Big) &= -\frac{1}{|\ub|^2} S Z_b^{l-1} Z_g^{k-l}\phi+\frac{1}{\ub}L S Z_b^{l-1} Z_g^{k-l}\phi\\
 &= \underbrace{\frac{1}{|\ub|^2} S Z_b^{l-1} Z_g^{k-l}\phi}_{A}+\underbrace{\frac{1}{\ub}L Z_b^{l-1} Z_g^{k-l+1}\phi}_{B}.
\end{align*}
For $A$, by rewriting $S$ as $\ub L + u\Lb$, we can use induction hypothesis for $(l, k-1)$; for $B$, we can use induction hypothesis for 
$(l-1, k)$.

Hence, we have obtained the desired estimates for $L Z_b^{l} Z_g^{k-l}\phi$.

\bigskip

Finally, to bound $\Lb Z_b^{l} Z_g^{k-l}\phi$, we use the equation
\begin{align*}
 -L\Lb Z_b^{l} Z_g^{k-l} \phi+ \slashed{\triangle} Z_b^{l} Z_g^{k-l}\phi =-\frac{1}{r}(L Z_b^{l} Z_g^{k-l}\phi-\Lb Z_b^{l} Z_g^{k-l}\phi) + \sum_{p+q\leq k} Q(\nabla Z^p\phi, \nabla Z^q \phi).
\end{align*}
We rewrite this as
\begin{align*}
 &\ \ -L\Lb Z_b^{l} Z_g^{k-l} \phi-\frac{1}{r} \Lb Z_b^{l} Z_g^{k-l}\phi + Q(\nabla \phi, \nabla Z_b^{l} Z_g^{k-l} \phi)\\
 &=-\slashed{\triangle} Z_b^{l} Z_g^{k-l}\phi-\frac{1}{r} L Z_b^{l-1} Z_g^{k-l}\phi + \sum_{p+q\leq k \atop
 p<k, q<k} Q(\nabla Z^p\phi, \nabla Z^q \phi).
\end{align*}
All the terms on the right-side have been controlled in previous steps. Therefore, it is straightforward to see that the right-hand side 
is bounded by $C(I_{n+1})|\ub|^{-2}\delta^\frac{1}{4}$. We now can mimic the proof for \eqref{a1} by defining $y=\ub\Lb Z_b^{l} Z_g^{k-l} \phi$, this leads to
\begin{equation*}
 |\Lb Z_b^{l} Z_g^{k-l}\phi(\ub,\delta,\theta)-\frac{C(1-\delta}{\ub}\Lb Z_b^{l} Z_g^{k-l}\phi(1-\delta,\delta,\theta)| \lesssim \delta^{\frac{1}{4}}|\ub|^{-1}C(I_{n+1}).
\end{equation*}
Taking into account of the vanishing property of $\Lb Z_b^{l} Z_g^{k-l}\phi$ on $S_{1-\delta,\delta}$, we complete the proof of Proposition \ref{impoved L infty estmiates}.

\begin{remark}
 For applications in the next section, we only need a slightly weakened (in decay) version of the estimates from Proposition \ref{impoved L infty estmiates}:
 \begin{equation}\label{impoved L infty estmiates for pratical uses}
\begin{split}
 \|\nablaslash Z_b^l Z_g^{k-l}\phi \|_{L^\infty(S_{\ub,\delta})}+ \|L Z_b^l Z_g^{k-l}& \phi\|_{L^{\infty}(S_{\ub,\delta})}\lesssim \delta^{\frac{1}{4}}|\ub|^{-\frac{3}{2}}C(I_{n+1}),\\
\|\Lb Z_b^l Z_g^{k-l} \phi\|_{L^\infty(S_{\ub,\delta})} &\lesssim \delta^{\frac{1}{4}}{|\ub|}^{-1} C(I_{n+1}),
\end{split}
\end{equation}
where $0 \leq l \leq k \leq n-3$.
\end{remark}

\section{Small data region}

In this section, we construct solutions in the entire small data region, i.e. region I. The approach is a modification of the classical approach with additional difficulties arising from the boundary  $C_\delta$.

\subsection{Klainerman-Sobolev inequality revisited}
We first introduce notations needed for the statement of the Klainerman-Sobolev inequality. We use $\Sigma_t$ to denote the constant time slices in the small data region, i.e. for a fixed $t \in (1,+\infty)$,
\begin{equation*}
 \Sigma_t := \big\{(x,t) \big| t-r \geq \delta \big\}.
\end{equation*}
This is a ball of radius $t-\delta$. We recall that we use $\Sigma_{1}$ to denote the entire $t=1$ hyperplane. Given a point $(t,x) \in \Sigma_t$ (assuming that $x \neq 0$), we use the $(t,B(t,x))$ to denote its corresponding boundary point, i.e. $(t,B(t,x))$ is the unique point on the boundary of $\Sigma_t$ (also on $C_\delta$)which is the intersection of the boundary of $\Sigma_t$ with the ray emanated from $(t,0)$ and passing from $(t,x)$. We now state the Klainerman-Sobolev inequality:
\begin{proposition}For all $f \in C^{\infty}(\mathbb{R}^{3+1})$, $t > 1$ and a point $(t,x)$ in the small data region, we have
\begin{equation}\label{KlaiermanSobolev}
|f(t,x)|\lesssim\frac{1}{(1+|u|)^{1/2}}|f(t,B(t,x))|+\frac{1}{(1+|\ub|)(1+|u|)^{1/2}}\sum_{Z \in \mathcal{Z}, k\leq3}\|Z^{k} f\|_{L^{2}(\Sigma_{t})}.
\end{equation}
\end{proposition}

We recall the following identities on $\mathbb{R}^{3+1}$:
\begin{equation}\label{formulas for partial}
\begin{split}
\partial_t &= \frac{1}{t-r} \big(\frac{t}{t+r} S - \sum_{i=1}^3\frac{x^i}{t+r} \Omega_{0i}\big), \\
\partial_i &= -\frac{1}{t-r} \big(\frac{x^i}{t+r} S -\frac{t}{t+r} \Omega_{0i} - \sum_{j=1}^3\frac{x^j}{t+r} \Omega_{ij}\big),\\
\partial_r &= \frac{1}{t-r} \big(-\frac{r}{t+r} S + \sum_{i=1}^3\frac{tx^i}{(t+r)r} \Omega_{0i}\big).
\end{split}
\end{equation}
Therefore, schematically, in terms of $Z_g \in \mathcal{Z}_g$, we write the above identities as
\begin{equation*}
\partial =\frac{1}{|t-r|}Z_g
\end{equation*}
Near light cone $C_0$. i.e. the hypersurface $t=r$, we can take $Z$ to be $\partial_i$ or $\partial_t$, therefore, schematically we have
\begin{equation*}
\partial =\left(1+\frac{1}{u}\right)Z.
\end{equation*}
We remark that this schematic expression means, for any function $f$, we have the following pointwise estimates:
\begin{equation*}
|\partial f|\lesssim \left(1+\frac{1}{u}\right)|Zf|.
\end{equation*}

We start the proof of \eqref{KlaiermanSobolev}.
Let $\chi$ be a non-negative smooth cut-off function on $\mathbb{R}_{\geq 0}$ so that $\chi$ is supported in $[0,\frac{1}{2}]$ and $\chi \equiv 1$ on $[0, \frac{1}{4}]$. We decompose $f(t,x)$ as
\begin{align*}
f(t,x) &= f_1(t,x) + f_2(t,x)\\
& = \chi(\frac{x}{t})f(t,x) + (1-\chi(\frac{x}{t}))f(t,x).
\end{align*}
Therefore, the function $f_1(t,x)$ is supported in region
\begin{equation*}
D_1 = \big\{(t,x) \mid 2r \leq t, \, t\geq 1\big\},
\end{equation*}
which is far away from the cone $C_{\delta}$; the function $\phi_2(t,x)$ is supported in region
\begin{equation*}
D_2 = \big\{(t,x) \mid t-r \geq \delta,\, 4r \geq t, \, t\geq 1\big\}
\end{equation*}
which is close to the cone $C_{\delta}$.

\bigskip

We first bound $f_1(t,x)$ in region $D_1$. In the rest of the subsection, we regard $t$ as a fixed large parameter. Let
\begin{equation*}
\widetilde{f_{1}}(x)=f_1(t,tx)=f(t,tx)\chi(x),
\end{equation*}
therefore, for a given positive integer $m$, we have
\begin{equation*}
\begin{split}
\|\partial^{m}\widetilde{f_{1}}\|^{2}_{L^{2}(\mathbb{R}^3)}&=\int_{\mathbb{R}^{3}}|\partial^{m}(f(t,tx)\chi(x))|^{2}dx \\
&\lesssim \sum_{j=0}^m \int_{\mathbb{R}^{3}} |\partial^{j}(f(t,tx))|^{2}|\nabla^{m-j}\chi|^{2}dx\\
&\lesssim \sum_{j=0}^m  \int_{\mathbb{R}^{3}} |t^{j}(\partial^{j}f)(t,tx)|^{2} dx.
\end{split}
\end{equation*}
Recall that we have $\partial = \frac{1}{1+|t-r|}{Z}$. In the region $D_1$, we have $t \geq 2 r$, hence $|t-r| \sim t$. Therefore, in $D_1$, we have
\begin{equation*}
|t\partial f|\lesssim |Zf|.
\end{equation*}
Thus, we have
\begin{align*}
\|\partial^{m}\widetilde{f_{1}}\|^{2}_{L^{2}(\mathbb{R}^3)} &\lesssim \sum_{j \leq m, Z \in \mathcal{Z}}\int_{\mathbb{R}^{3}}|Z^{j}f|^{2}(t,tx)dx \\
&= t^{-3} \sum_{j \leq m, Z \in \mathcal{Z}} \int_{\mathbb{R}^{3}}|Z^{j}f(t,y)|^{2}dy.
\end{align*}
Therefore, according to the classical Sobolev inequality on $\mathbb{R}^3$, we obtain
\begin{equation}\label{estimates in region C_1}
\|f_1\|_{L^\infty(\Sigma_t)} = \|\widetilde{f_1}\|_{L^\infty(\Sigma_t)} \lesssim \frac{1}{t^{\frac{3}{2}}} \sum_{k \leq 2, Z \in \mathcal{Z}}\|Z^{k} f(t,\cdot)\|_{L^{2}(\Sigma_t)}.
\end{equation}
\bigskip

We turn to the estimates on 
$f(t,x)$ in the region $D_2$. On the hyperplane $\Sigma_t$, we draw a line from the origin and the point $(x,t)$. When a point moves along the radial direction on this line, it hits the characteristic boundary of $C_{\delta}$ at one point $(t,B(t,x))$. By integrating $\partial_r \left((1+|t-r|) f^{2}(t,x)\right)$  from $(t,B(t,x))$ to $(t,x)$, we obtain
\begin{equation*}
\begin{split}
(1+|t-r|)f^{2}(t,x)&=(1+\delta)f^{2}(t,B(t,x))+\int_{r}^{t-\delta}\partial_{r}\big( (1+|t-r|)f^{2}(t,x)\big)dr\\
&=(1+\delta)f^{2}(t,B(t,x))+\int_{r}^{t-\delta}-f^{2}(t,x) + 2 (1+|t-r|) f(t,x)\partial_{r}f(t,x)dr
\end{split}
\end{equation*}
For the integrand in the last line, we apply the classical Sobolev inequalities on spheres $S_{t,r}$ (the sphere of radius $r$ on $\Sigma_t$). Therefore, we obtain
\begin{align*}
((1+|t-r|)f^{2}(t,x)&\lesssim f^{2}(t,B(t,x)) + \int_{r}^{t-\delta} \frac{1}{r^2}\sum_{|\alpha|\leq2}\|\Omega^{\alpha}f\|_{L^{2}(S_{t,r})}^{2}dr\\
&\ \ +\int_{r}^{t-\delta} \frac{(1+|t-r|)}{r^2}\sum_{|\alpha|, |\beta| \leq 2}\|\Omega^{\alpha}f\|_{L^2(S_{t,r})}\|\Omega^{\beta}\partial_{r}f\|_{L^2(S_{t,r})}dr.
\end{align*}
Since $|t-r|\partial_{r}\lesssim Z$, we have
\begin{align}\label{estimates in region C_2}
\begin{split}
((1+|t-r|)f^{2}(t,x)&\lesssim f^{2}(t,B(t,x)) +  \int_{r}^{t-\delta}  \frac{1}{r^2}\sum_{k\leq 3, Z \in \mathcal{Z}}\|Z^{k}f\|^2_{L^2(S_{t,r})}dr\\
&=f^{2}(t,B(t,x)) + \frac{1}{r^2}  \sum_{k|\leq 3, Z \in \mathcal{Z}}\|Z^{k} f(t,\cdot)\|^2_{L^2(\Sigma_t)}.
\end{split}
\end{align}

The estimates \eqref{estimates in region C_1} and \eqref{estimates in region C_2} together give the desired estimates \eqref{KlaiermanSobolev} and we complete the proof.

\subsection{A priori energy estimates}

For a $t \in (1,+\infty)$, we still use $\Sigma_t = \big\{(x,t) \big| t-r \geq \delta \big\}.$ to denote the constant time slices in the small data region.  For $k \in \mathbb{Z}_{\geq 0}$ and $t>1$, we introduce the following energy norms:
\begin{equation}
\begin{split}
\widetilde{E}_k(t) &= \big(\sum_{Z \in \mathcal{Z}} \int_{\Sigma_t} |\partial_t Z^k \phi|^2+ \sum_{j=1}^3 |\partial_j Z^k \phi|^2 dx\big)^\frac{1}{2},\\
\widetilde{E}_{\leq k}(t) &= \big(\sum_{0\leq j \leq  k} \widetilde{E}_{j} (t)^2\big)^\frac{1}{2}.
\end{split}
\end{equation}

We use $D_{t,\delta}$ to denote the space-time region bounded by $\Sigma_t$, $\Sigma_0$ and $C_\delta$. This region is obviously foliated by the constant time foliation $\{\Sigma_\tau \mid \tau \in [1,t]\}$ and this foliation is one of the foliations we use to derive energy estimates. The second foliation is the null foliation of outgoing null cones $\{C_u \mid u\in [\delta,t/2]\}$. This foliation is depicted as follows:

\includegraphics[width = 5 in]{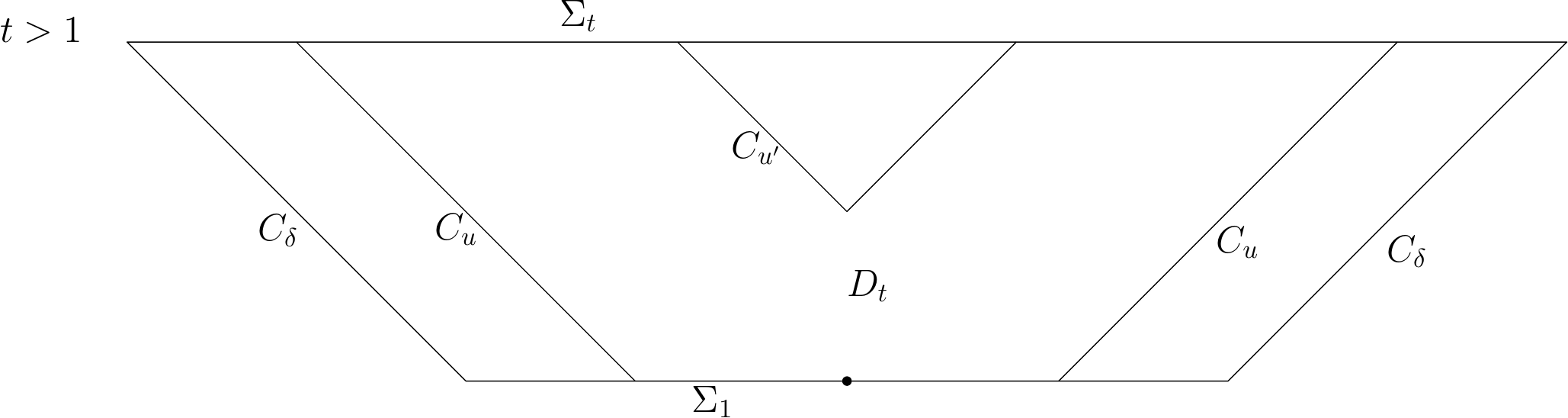}

Whenever there is no confusion, we still use $C_u$ to denote $C_u \cap D_t$. We use $D_{t,u}$ to denote the space-time region bounded by $\Sigma_t$, $\Sigma_1$ and $C_u$. This is a truncated solid light cone in $\mathbb{R}^{3+1}$. We use $\Sigma_{1,u}$ and $\Sigma_{t,u}$ to denote its bottom and top respectively. We remark that the bottom can be a single point.

Recall that (assuming that the solution $\phi$ exists up to time $t$), for $k\geq 0$ and $|\alpha|=k$, we have
\begin{equation*}
\Box Z^{k} \phi = \sum_{p+q \leq k}Q(\nabla Z^p \phi, \nabla Z^q \phi).
\end{equation*}
We multiply both sides by $\partial_t Z^k \phi$ and we then integrate over $D_{t,u}$. This leads to the following energy identity:
\begin{align*}
\int_{\Sigma_{t,u}}|\partial_t Z^k \phi|^2+ \sum_{j=1}^3 |\partial_j Z^k\phi|^2&= \int_{\Sigma_{1,u}}|\partial_t Z^k\phi|^2+ \sum_{j=1}^3 |\partial_j Z^k \phi|^2+\int_{C_u}|L Z^k\phi|^2 +|\nablaslash Z^k\phi|^2 \\
&\ \ + \sum_{p+q \leq k} \doubleint_{D_{t,u}} Q(\nabla Z^p \phi, \nabla Z^q \phi)\partial_t Z^k \phi.
\end{align*}
Recall that we use $\partial \in \{L,\Lb, \nablaslash\}$ to denote a generic derivative and use $\partial_g \in \{L,\nablaslash\}$ to denote a good derivative. Therefore, by using $|\partial Z^k\phi|^2$  as a shorthand notation for $ |\partial_t Z^k\phi|^2+ \sum_{j=1}^3 |\partial_j Z^k \phi|^2$ and using $|\partial_g Z^k\phi|^2$  as a shorthand notation for $ |L Z^k\phi|^2+ |\nablaslash Z^k \phi|^2$, we have
\begin{equation*}
\int_{\Sigma_{t,u}}|\partial Z^k \phi|^2 = \int_{\Sigma_{1,u}}|\partial Z^k\phi|^{2}+\int_{C_u}|\partial_g Z^k\phi|^2 + \sum_{p+q \leq k} \doubleint_{D_{t,u}} Q(\nabla Z^p \phi, \nabla Z^q \phi)\partial_t Z^k \phi.
\end{equation*}
In applications, since the data prescribed on 
$\Sigma_{1,u}$ are trivial, we have
\begin{equation}\label{Fundamental Energy Estimate}
\int_{\Sigma_{t,u}}|\partial Z^k \phi|^2 =\int_{C_u}|\partial_g Z^k\phi|^2 + \sum_{p+q \leq k} \doubleint_{D_{t,u}} Q(\nabla Z^p \phi, \nabla Z^q \phi)\partial_t Z^k \phi.
\end{equation}
\bigskip

Before we state the main estimates of the section, we first compute the energy flux $\int_{C_\delta}|\partial_g Z^k\phi|^2$ through the outermost cone $C_\delta$. According to \eqref{impoved L infty estmiates}, for $k\leq n-2$, we have $|\partial_g Z^k \phi|\lesssim \delta^\frac{1}{4}|\ub|^{-\frac{3}{2}-\frac{\alpha}{4}}C(I_{n+1})$, therefore,
\begin{equation}\label{flux through C delta}
\int_{C_\delta}|\partial_g Z^k\phi|^2 \lesssim \delta^\frac{1}{2}C(I_{n+1}),
\end{equation}
where we still use $C(I_{n+1})$ to denote $C(I_{n+1})^2$.

\begin{proposition}
Under the same assumptions as in the previous section, for sufficiently small $\delta$, there exists a unique global future in time solution $\phi$ of \eqref{Main Equation} on the small data region, so that together with the solution constructed in the short pulse region, we have a unique future in time solution $\phi$. Moreover, this solution $\phi$ on the small dat region enjoys the following energy estimates:
\begin{equation}\label{energy bound for solution}
\widetilde{E}_{\leq n}(t) \lesssim  \delta^\frac{1}{4}C(I_{n+1}),
\end{equation}
for all $t >1$,
\end{proposition}
\begin{remark}
The existence of solutions in the small data region follows from the \emph{a priori} estimate \eqref{energy bound for solution}. Since the argument is routine, we will not pursuit this point here.
\end{remark}
We use a bootstrap argument to prove the proposition. We assume that the solution exists up to time $t$ and for all $1\leq t' \leq t$, we have
\begin{equation}\label{bootstrap assumption for small data region}
\widetilde{E}_{\leq 7}(t') \lesssim  M \delta^\frac{1}{4}.
\end{equation}
It suffices to show that we can indeed choose $M$ so that it depends only on $I_n$

\bigskip
We first point out that we can derive $L^\infty$ bound on $\partial  Z^p \phi$ for $p\leq 4$. According to Klainerman-Sobolev inequality, we have
\begin{align*}
|\partial Z^p \phi(\tau,x)|&\lesssim \frac{1}{(1+|u|)^{1/2}}|\partial Z^p \phi(\tau,B(\tau,x))|+\frac{1}{(1+|\ub|)(1+|u|)^{1/2}}\sum_{Z \in \mathcal{Z}, l\leq3}\|Z^{l} \partial Z \phi\|_{L^{2}(\Sigma_{\tau})}\\
&\lesssim \frac{C(I_{n+1})}{(1+\ub)(1+u)^{1/2}}\delta^\frac{1}{4}+\frac{M}{(1+\ub)(1+|u|)^{1/2}}\delta^\frac{1}{4}.
\end{align*}
In particular, based on \eqref{formulas for partial}, it is well known that for good derivatives $\partial_g$, we have
\begin{align*}
|\partial_g Z^p \phi(\tau,x)|\lesssim \frac{M}{t^\frac{3}{2}}\delta^\frac{1}{4}.
\end{align*}

\bigskip

For all $u \geq \delta$, according to  \eqref{Fundamental Energy Estimate}, we have	
\begin{align*}
\int_{C_u}|\partial_g Z^k\phi|^2 &\leq \int_{\Sigma_{t,u}}|\partial Z^k \phi|^2 + \sum_{p+q \leq k} \doubleint_{D_{t,u}} \big|Q(\nabla Z^p \phi, \nabla Z^q \phi)\partial_t Z^k \phi\big|\\
&\leq \int_{\Sigma_{t}}|\partial Z^k \phi|^2 + \sum_{p+q \leq k} \doubleint_{D_{t,\delta}} \big|Q(\nabla Z^p \phi, \nabla Z^q \phi)\partial_t Z^k \phi\big|.
\end{align*}
For the last step, we have enlarged the domain for integration. Therefore, according to the foliation $\{u\in[\delta, \frac{t}{2}]\mid C_u\}$, for the given constant $\varepsilon_0 \in (0,\frac{1}{2})$, we have
\begin{align*}
 \doubleint_{D_{t,\delta}}\frac{|\partial_g Z^k\phi|^2 }{(1+|u|)^{1+\varepsilon_0}} &= \int_{\delta}^{t/2}\frac{1}{(1+|u|)^{1+\varepsilon_0}} \big(\int_{C_{u}} |\partial_g Z^k\phi|^2\big)du\\
 &\leq\int_{\delta}^{t/2}\frac{1}{(1+|u|)^{1+\varepsilon_0}} \big(\int_{\Sigma_{t}}|\partial Z^k \phi|^2 + \sum_{p+q \leq k} \doubleint_{D_{t,\delta}} \big|Q(\nabla Z^p \phi, \nabla Z^q \phi)\partial_t Z^k \phi\big|\big)du'
\end{align*}
Since the quantity inside the parenthesis is independent of $u'$, we obtain
\begin{equation}\label{aa}
 \doubleint_{D_{t,\delta}}\frac{|\partial_g Z^k\phi|^2 }{(1+|u|)^{1+\varepsilon_0}} \lesssim \int_{\Sigma_{t}}|\partial Z^k \phi|^2 + \sum_{p+q \leq k} \doubleint_{D_{t,\delta}} \big|Q(\nabla Z^p \phi, \nabla Z^q \phi)\partial_t Z^k \phi\big|.
\end{equation}
We take $u=\delta$ in \eqref{Fundamental Energy Estimate}. In view of \eqref{flux through C delta}, we obtain immediately that
\begin{equation}\label{bb}
\begin{split}
 \int_{\Sigma_{t}}|\partial Z^k \phi|^2 &\leq \int_{C_\delta}|\partial_g Z^k\phi|^2 + \sum_{p+q \leq k} \doubleint_{D_{t,\delta}} |Q(\nabla Z^p \phi, \nabla Z^q \phi)||\partial_t Z^k \phi|\\
 &\lesssim \delta^\frac{1}{2}C(I_{n+1}) +  \sum_{p+q \leq k} \doubleint_{D_{t,\delta}} |Q(\nabla Z^p \phi, \nabla Z^q \phi)||\partial_t Z^k \phi|.
\end{split}
\end{equation}
Together with \eqref{aa}, we have
\begin{equation}\label{cc}
 \doubleint_{D_{t,\delta}}\frac{|\partial_g Z^k\phi|^2 }{(1+|u|)^{1+\varepsilon_0}} \lesssim \delta^\frac{1}{2}C(I_{n+1}) +  \sum_{p+q \leq k} \doubleint_{D_{t,\delta}} |Q(\nabla Z^p \phi, \nabla Z^q \phi)||\partial_t Z^k \phi|.
 \end{equation}
In view of \eqref{bb}, we arrive at the following energy estimates:
\begin{equation*}
 \int_{\Sigma_{t}}|\partial Z^k \phi|^2 + \doubleint_{D_{t,\delta}}\frac{|\partial_g Z^k\phi|^2 }{(1+|u|)^{1+\varepsilon_0}} \lesssim \delta^\frac{1}{2}C(I_{n+1}) +  \sum_{p+q \leq k} \doubleint_{D_{t,\delta}} |Q(\nabla Z^p \phi, \nabla Z^q \phi)||\partial_t Z^k \phi|.
 \end{equation*}
By summing over $k$, we finally obtain that
\begin{equation}\label{modified energy identity}
 \widetilde{E}_{\leq k}(t) + \sum_{l\leq k,\atop Z \in \mathcal{Z}}\doubleint_{D_{t,\delta}}\frac{|\partial_g Z^l\phi|^2 }{(1+|u|)^{1+\varepsilon_0}} \lesssim \delta^\frac{1}{2}C(I_{n+1}) +  \sum_{l \leq k, p+q \leq l, \atop Z \in \mathcal{Z}}\!\!\! \doubleint_{D_{t,\delta}} |Q(\nabla Z^p \phi, \nabla Z^q \phi)||\partial_t Z^l\phi|.
 \end{equation}
Since we have the energy term $\widetilde{E}_{\leq k}(t)$ on the left-hand side, to complete the bootstrap argument, it suffices to control the second term on the right-hand side. According to the structure of null forms, this term is bounded by
\begin{align*}
 \sum_{l \leq k, p+q \leq l, \atop Z \in \mathcal{Z}}\!\!\! \doubleint_{D_{t,\delta}} |\partial_g Z^p \phi||\partial Z^q \phi||\partial_t Z^l\phi|.
\end{align*}
According to whether $p < q$ or $p \geq q$, we break this term into two pieces (we replace $\partial_t$ by $\partial$):
\begin{align*}
 S_1 + S_2 = \sum_{l \leq k, p+q \leq l, \atop p < q, Z \in \mathcal{Z}}\!\!\! \doubleint_{D_{t,\delta}} |\partial_g Z^p \phi||\partial Z^q \phi||\partial Z^l\phi|+
 \sum_{l \leq k, p+q \leq l, \atop p \geq q, Z \in \mathcal{Z}}\!\!\! \doubleint_{D_{t,\delta}} |\partial_g Z^p \phi||\partial Z^q \phi||\partial Z^l\phi|.
\end{align*}

For $S_1$, since $p<q$, we have  $k-p \geq \lfloor \frac{1}{2}k \rfloor \geq 3$. Here $\lfloor \frac{1}{2}k \rfloor$ denotes the largest integer less or equal to $\frac{1}{2}$. We can apply the $L^\infty$ estimates for good derivatives $\partial_g Z^p \phi$. Therefore,
\begin{align*}
 S_1 &\lesssim \sum_{l \leq k, p+q \leq l, \atop p < q, Z \in \mathcal{Z}}\int_{1}^t \frac{M}{\tau^\frac{3}{2}} \delta^\frac{1}{4}\|\partial Z^q \phi\|_{L^2(\Sigma_\tau)}\|\partial Z^l\phi\|_{L^2(\Sigma_\tau)}d\tau\\
 &\lesssim M^3 \delta^\frac{3}{4}.
\end{align*}

For $S_2$, we apply Klainerman-Sobolev to $|\partial Z^q \phi|$ and we obtain
\begin{align*}
 S_2 &\lesssim \sum_{l \leq k, p+q \leq l, \atop p \geq q, Z \in \mathcal{Z}}\!\!\! \doubleint_{D_{t,\delta}}\frac{M}{t(1+|u|)^\frac{1}{2}} \delta^\frac{1}{4}|\partial_g Z^p \phi||\partial Z^l\phi|\\
 &\lesssim \epsilon \sum_{p \geq k, Z \in \mathcal{Z}}\doubleint_{D_{t,\delta}}\frac{|\partial_g Z^p\phi|^2 }{(1+|u|)^{1+\varepsilon_0}}+\frac{1}{\epsilon}\sum_{l \leq k, p+q \leq l, \atop p \geq q, Z \in \mathcal{Z}}\!\!\! \doubleint_{D_{t,\delta}}\frac{M^2(1+|u|)}{t^2} \delta^\frac{1}{2}|\partial Z^l\phi|^2\\
 &\lesssim  \epsilon \sum_{p \geq k, Z \in \mathcal{Z}}\doubleint_{D_{t,\delta}}\frac{|\partial_g Z^p\phi|^2 }{(1+|u|)^{1+\varepsilon_0}} + \frac{1}{\epsilon} M^4 \delta,
\end{align*}
where the constant $\epsilon$ will be determined later on.

Back to \eqref{modified energy identity}, the estimates on $S_1$ and $S_2$ yield
\begin{equation*}
 \widetilde{E}_{\leq k}(t) + \sum_{l\leq k,\atop Z \in \mathcal{Z}}\doubleint_{D_{t,\delta}}\frac{|\partial_g Z^l\phi|^2 }{(1+|u|)^{1+\varepsilon_0}} \lesssim \delta^\frac{1}{2}C(I_{n+1}) +M^3 \delta^\frac{3}{4}+ \epsilon \sum_{p \geq 7, Z \in \mathcal{Z}}\doubleint_{D_{t,\delta}}\frac{|\partial_g Z^p\phi|^2 }{(1+|u|)^{1+\varepsilon_0}} + \frac{1}{\epsilon} M^4 \delta.
 \end{equation*}
By choosing a suitable small constant $\epsilon$, we can remove the integral term on the right-hand side and obtain
\begin{equation*}
 \widetilde{E}_{\leq k}(t) + \sum_{l\leq k,\atop Z \in \mathcal{Z}}\doubleint_{D_{t,\delta}}\frac{|\partial_g Z^l\phi|^2 }{(1+|u|)^{1+\varepsilon_0}} \lesssim \delta^\frac{1}{2}C(I_{n+1}) +M^3 \delta^\frac{3}{4} + \frac{1}{\epsilon} M^4 \delta.
 \end{equation*}
Hence,
 \begin{equation*}
 \widetilde{E}_{\leq k}(t) \lesssim \delta^\frac{1}{2}C(I_{n+1}) +M^3 \delta^\frac{3}{4} + \frac{1}{\epsilon} M^4 \delta.
 \end{equation*}
We then can choose a sufficiently small $\delta$ and this completes the bootstrap argument.

\section*{Acknowledgment}
S.Miao is supported by NSF grant DMS-1253149 to The University of Michigan; L.Pei is supported by the PhD research fellowship of Norwegian University of Science and Technology; P.Yu is supported by NSFC 11101235 and NSFC 11271219. S.Miao would like to thank Sohrab Shahshahani for helpful comments on a previous version of this manuscript. P.Yu would like to thank Sergiu Klainerman for the communication on the relaxation of the propagation estimates.

\end{document}